\documentclass[smallextended,referee,envcountsect]{svjour3} 
\smartqed 
\journalname{JOTA}

\usepackage{algorithm}
\usepackage{algorithmic}
\usepackage{subcaption}
\usepackage{float}
\usepackage{graphicx}
\usepackage{xcolor}
\usepackage{hyperref}       
\usepackage{url}

\usepackage{amsmath, amssymb}
\newtheorem{assumption}{Assumption}
\newcommand*{\R}{{\mathbb R}}
\newcommand*{\E}{{\mathbb E}}
\DeclareMathOperator*{\argmin}{argmin}
\DeclareMathOperator*{\argmax}{argmax}

\newcommand*{\Z}{{\mathcal Z}}
\newcommand*{\PP}{{\mathbb P}}

\def\ag#1{{\color{black}#1}} 
\def\aai#1{{\color{black}#1}} 
\def\ai#1{{\color{black}#1}} 
\def\pd#1{{\color{black}#1}} 
\def\att#1{{\color{black}#1}} 

\def\rev#1{{\color{black}#1}} 

 \journalname{JOTA}
\begin{document}

\title{Oracle Complexity Separation in Convex Optimization
}


\author{Anastasiya Ivanova \and
Pavel Dvurechensky \and
Evgeniya Vorontsova \and
Dmitry Pasechnyuk \and
Alexander Gasnikov \and
Darina Dvinskikh \and 
Alexander Tyurin
}


\institute{A. Ivanova \at
              corresponding author,
              National Research University Higher School of Economics, Moscow, Russian Federation; \\
              Grenoble Alpes University, Grenoble, France \\
              \email{asivanova@hse.ru} 
           \and
Pavel Dvurechensky 
\at 
Weierstrass Institute for Applied Analysis and Stochastics, Berlin, Germany \\
\email{pavel.dvurechensky@wias-berlin.de}
\and
Evgeniya Vorontsova
\at Institute of Information and Communication Technologies, Electronics and Applied Math., Catholic University of Louvain,
Belgium; \\
Moscow Institute of Physics and Technology, Moscow, Russian Federation
\email{vorontsovaea@gmail.com}
\and
Dmitry Pasechnyuk 
\at Moscow Institute of Physics and Technology, Moscow, Russian Federation\\
ISP RAS Research Center for Trusted Artificial Intelligence, Moscow, Russian Federation\\
\email{pasechnyuk2004@gmail.com}
\and
Alexander Gasnikov 
\at Moscow Institute of Physics and Technology, Moscow, Russian Federation;\\ 
ISP RAS Research Center for Trusted Artificial Intelligence, Moscow, Russian Federation\\  
Institute for Information Transmission Problems, Moscow, Russian Federation\\
\email{gasnikov@yandex.ru}
\and
Darina Dvinskikh 
\at 
ISP RAS Research Center for Trusted Artificial Intelligence, Moscow, Russian Federation\\ 
Moscow Institute of Physics and Technology, Moscow, Russian Federation \\
\email{dviny.d@yandex.ru}
\and 
Alexander Tyurin
\at
National Research University Higher School of Economics, Moscow, Russian Federation \\
\email{alexandertiurin@gmail.com}
}

\date{Received: date / Accepted: date}

\maketitle

\begin{abstract}
Many convex optimization problems have structured objective functions written as a sum of functions with different types of oracles (e.g., full gradient, coordinate derivative, stochastic gradient) and different arithmetic operations complexity of these oracles. In the strongly convex case, these functions also have different condition numbers that eventually define the iteration complexity of first-order methods and the number of oracle calls required to achieve a given accuracy. Motivated by the desire to call more expensive oracles fewer times, we consider the problem of minimizing a sum of two functions and propose a generic algorithmic framework to separate oracle complexities for each function. The latter means that the oracle for each function is called the number of times that coincide with the oracle complexity for the case when the second function is absent.
Our general accelerated framework covers the setting of (strongly) convex objectives, the setting when both parts are given through full coordinate oracle, as well as when one of them is given by coordinate derivative oracle or has the finite-sum structure and is available through stochastic gradient oracle. In the latter two cases, we obtain accelerated random coordinate descent and accelerated variance reduced methods with oracle complexity separation.

\end{abstract}

\section{Introduction}
\label{intro}
The complexity of an optimization problem usually depends on the objective parameters, such as the Lipschitz constant of the gradient and the strong convexity parameter. 
In this paper, we are mainly motivated by Machine Learning applications, where the objective is constructed from many building blocks, e.g., individual loss for an example in the dataset or different regularizers in supervised machine learning. A similar optimization problem structure is also used in imaging problems.
Standard theoretical results about optimization algorithms for such problems provide \textit{iteration} complexity, namely, the number of iterations to achieve a given accuracy. Unlike these results, in this paper, we address the question of \textit{oracle} complexity, focusing on the number of oracle calls. Moreover, the goal is to study separately what number of oracle calls \textit{for each} building block of the objective is sufficient to obtain the required accuracy. Indeed, typically in the empirical loss minimization, the finite-sum part of the objective is much more computationally expensive than the regularizer, which motivates the usage of a randomized oracle for the finite-sum part and a proximal oracle for the regularizer. Further, some components in the finite-sum part may be more computationally expensive than others, and it is desirable to call the gradient oracle of the former less frequently than the gradient oracle of the latter. Moreover, some of the building blocks of the objective may be available through their gradient, while for the other blocks, only the value of the objective may be available. In this case, one would prefer to call the gradient oracle for the former less frequently than the zeroth-order oracle for the latter. 
Even in the first-order oracle model, the evaluation of the first-order oracle for one term in the sum can be more expensive than for another. 
Thus, if we manage to estimate separately for each building block of the objective the number of oracle calls that is sufficient to achieve some desired accuracy, and the oracle evaluation for each block has a different cost, we can hope that more expensive oracles will be called less frequently, and the total working time of the algorithm will be smaller.
This paper proposes a generic framework that allows achieving such oracle complexity separation for accelerated first-order methods in different contexts. 

\aai{We are aware of only a few works in this direction, which, to our knowledge, was initiated in \cite{lan2016gradient}, where complexity separation is proposed for a sum of a smooth and a non-smooth functions.}
In \cite{lan2016accelerated} the authors propose an accelerated algorithm to separate complexities in the deterministic full gradient setting both with convex and strongly convex objectives. Unlike them, we construct a generic framework that allows one to work with other types of oracles such as stochastic gradient in finite-sum optimization or coordinate derivative, which allows one to construct accelerated randomized methods.
In \cite{lan2016conditional}, the authors obtain separate estimates for the number of linear minimization oracle calls for a compact set and the number of gradient evaluations for a variant of the Frank-Wolfe method.

Another line of related works is on generic acceleration frameworks \cite{monteiro2013accelerated,lin2015universal,lin2018catalyst}. In particular, \cite{lin2015universal,lin2018catalyst} propose a generic scheme to accelerate non-accelerated methods, including randomized methods. Despite proposing generic frameworks, these papers do not consider the question of oracle complexity separation as we do. Our framework is inspired by the accelerated algorithm in \cite{monteiro2013accelerated}, which we generalize in several aspects.

The paper is organized as follows. In the next Section \ref{S:contrib}, we list main assumptions, give a detailed description of our contributions and compare them with the existing literature. In Section \ref{sec:1}, we describe our general framework for composite optimization and state the main result. Section \ref{appl} is devoted to three particular applications of our general framework: composite gradient methods, composite random coordinate descent, composite stochastic variance reduced methods. In Section \ref{exp}, we illustrate our framework by numerical experiments on Kernel Support Vector Machine problem and log-density estimation using Bayesian approach. The proof of the main result is provided in Section \ref{App:proof}.

\section{Problem Statement and Contributions}
\label{S:contrib}

\aai{The main optimization problem that we consider in this paper is
\begin{equation}
\label{prob_st}
    \min_{x\in\R^n} \{f(x):= h(x) + g(x)\},
\end{equation}
where $h$ and $g$ are convex and continuously differentiable functions.
We also make the following assumptions about the functions $f$, $h$ and $g$ denoting by $\|\cdot\|_2$ the standard Euclidean norm in $\R^n$.
\begin{assumption}
    The function $f(\cdot)$ is $L_f$-smooth, i.e. has Lipschitz continuous gradient with  constant $L_f$  w.r.t. $\|\cdot\|_2$, and is $\mu$-strongly convex w.r.t. $\|\cdot\|_2$ with $\mu \geq 0$.
\end{assumption}
\begin{assumption}
    The function $h(\cdot)$ is $L_h$-smooth, i.e. has Lipschitz continuous gradient with constant $L_h$ w.r.t. $\|\cdot\|_2$, and there is an oracle $O_h$ that in one call produces the gradient $\nabla h(\cdot).$
\end{assumption}
\begin{assumption}
\label{Asmpt:g}
   The function $g(\cdot)$ is $L_g$-smooth, i.e. has Lipschitz continuous gradient with constant $L_g$ w.r.t. $\|\cdot\|_2$, and there is a basic oracle $O_g$ that in $\kappa_{g}$ calls produces the gradient $\nabla g(\cdot)$.
\end{assumption}
}
Let us comment on the last assumption by giving three examples of the basic oracle. The most simple example is when $O_g$ returns $\nabla g(\cdot)$ leading to $\kappa_{g}=1$. The second example is when $O_g$ returns a coordinate derivative $\nabla_i g(\cdot)$ leading to $\kappa_{g}=n$ since to find the full gradient one needs to find coordinate derivatives for each of $n$ coordinates. Finally, if $g(x)=\tfrac{1}{m} \sum_{k=1}^{m} g_k(x)$ and $O_g$ returns  $\nabla g_k(\cdot)$, we have that $\kappa_{g}=m$.

Our framework is based on a series of inner-outer loops, where in each loop a problem with some specific structure is solved by some algorithm. In order to make our framework flexible and being able to work with different types of basic oracles $O_g$, we introduce below a generic assumption on an algorithm that is used in the innermost loop. 
\rev{In this loop we need to solve a problem that has the following form:
\begin{equation}
    \label{eq:innermost_problem}
    \min_{v \in \R^n} \left\{\varphi(v) = \left\langle \beta, v \right\rangle + \tfrac{\alpha}{2}\|v\|_2^2 + g(v)\right\},
\end{equation}
where $\beta \in \R^n$ is a given vector, $\alpha >0$ is a given number. We denote by $v^*$ the solution of this problem.
Strongly convex problem \eqref{eq:innermost_problem} is solved by an algorithm $\mathcal{M}_{inn}$ starting from a point $v^0$ and using $N$ calls to the  basic oracle $O_g$, which is denoted by $\mathcal{M}_{inn}(\varphi(\cdot),v^0,N)$. This algorithm is allowed to be randomized and output a random point $\hat{v}$. We use $\E$ to denote the expectation w.r.t. all the randomness introduced by the algorithm. The following assumption is made about the complexity of $\mathcal{M}_{inn}$ applied to solve problem \eqref{eq:innermost_problem}.  
}
\begin{assumption}
\label{ass_inn_method} 
\rev{There exist a parameter $\tau_{g} > \sqrt{\alpha}$ that depends on the function $g(\cdot)$ and the method $\mathcal{M}_{inn}$,  but is independent of $\alpha$, a constant $C>0$ and a numerical constant $C_0>0$ such that, for any accuracy $\varepsilon>0$, $\mathcal{M}_{inn}(\varphi(\cdot),v^0, N(\tilde{\varepsilon}))$ outputs a (possibly random) point $\hat{v}$ such that $\E \left(\varphi(\hat{v}) - \varphi(v^{\ast}) \right) \leq \tilde{\varepsilon}$, in $N(\tilde{\varepsilon}) = \tfrac{C_0\tau_{g}}{\sqrt{\alpha}}\ln\tfrac{C \|v^0 - v^{\ast}\|_2^2}{\tilde{\varepsilon}} $  calls to the basic oracle $O_g$.
}
\end{assumption}
\rev{This assumption may be illustrated by three examples. For that we consider $\left\langle \beta, v \right\rangle + \tfrac{\alpha}{2}\|v\|_2^2$ in \eqref{eq:innermost_problem} as a composite term.
Firstly, we can  apply Accelerated Gradient Method for Composite Optimization from~\cite{nesterov2013gradient} as $\mathcal{M}_{inn}$, which means that Assumption~\ref{ass_inn_method} holds with $\tau_g = \sqrt{L_{g}}$ since the smooth part $g$ of the objective is $L_g$-smooth and the whole objective $\varphi$ is $\alpha$-strongly convex. Note that in this case the point $\hat{v}$ is deterministic. Secondly, we can apply Accelerated Proximal Random Coordinate Method from~\cite{fercoq2015accelerated,allen2016even} as $\mathcal{M}_{inn}$, which means that Assumption~\ref{ass_inn_method} holds with $\tau_g =  n \sqrt{\overline{L}_{g}}$, where $\sqrt{\overline{L}_{g}} = \tfrac{1}{n}  \sum_{i=1}^n \sqrt{\beta_i}$, and $\beta_i$, $i=1,\ldots,n$ are Lipschitz constants for the coordinate derivatives of $g$. Finally, if $g(x)=\tfrac{1}{m} \sum_{k=1}^{m} g_k(x)$  we can apply Accelerated Stochastic Variance Reduced Algorithms, e.g. Katyusha~\cite{allen2017katyusha} or Varyag \cite{lan2019unified} as the inner method $\mathcal{M}_{inn}$. 
For these methods the number of oracle calls to solve the problem \eqref{eq:innermost_problem} is $\tilde{O}\left(m + \sqrt{\tfrac{m \hat{L}_g}{\alpha}} \right)$, where $\hat{L}_g = \max\limits_{k}L_{g_k}$ and Assumption~\ref{ass_inn_method} holds with $\tau_g =  \sqrt{m \hat{L}_g}$ if we consider ill-conditioned problems with $m\leq \sqrt{\tfrac{m \hat{L}_g}{\alpha}}$.
}

\rev{
Based on the above assumptions, 
we propose a general framework for solving  problem \eqref{prob_st} that, using some accelerated method to solve problems of the form \eqref{eq:innermost_problem},} provides an accelerated algorithm with separated oracle complexities for the parts $h$ and $g$. In particular, when $f$ is $\mu$-strongly convex, as a corollary of our general framework, we propose the following algorithms. 
\begin{enumerate}
    \item (Full gradient setting.) Under the assumption that $h$ is $L_h$-smooth and $g$ is $L_g$-smooth with $L_g \ge L_h$, we propose an Accelerated Gradient Method that obtains $\hat{x}$  such that $f(\hat{x}) - f(x^{\ast}) \leq \varepsilon$ in $O\left(\sqrt{\tfrac{L_h}{\mu }}\ln\left( \tfrac{1}{\varepsilon} \right) \right)$ gradient evaluations for $h(\cdot)$ and $ O \left(\sqrt{\tfrac{L_g}{\mu}} \ln\left( \tfrac{1}{\varepsilon} \right)  \right) $ gradient evaluations for $g(\cdot)$.
    \item (Random coordinate descent setting.) Under the assumption that $h$ is $L_h$-smooth and, for $i=1,...n$, $ \partial g(x) / \partial x_i$ is $\beta_i$-Lipschitz with $\sqrt{\overline{L}_{g}} = \tfrac{1}{n}  \sum_{i=1}^n \sqrt{\beta_i}$ and $\overline{L}_{g} \ge L_h$, we propose an Accelerated Random Coordinate Descent Method that obtains $\hat{x}$  such that $f(\hat{x}) - f(x^{\ast}) \leq \varepsilon$ with high probability in  $O\left(\sqrt{\tfrac{L_h}{\mu }} \ln\left( \tfrac{1}{\varepsilon} \right) \right)$ gradient evaluations for $h(\cdot)$ and $ O \left(n \sqrt{\tfrac{\overline{L}_{g}}{\mu}} \ln\left( \tfrac{1}{\varepsilon} \right)  \right) $ partial derivatives evaluations for $g(\cdot)$.
    \item (Variance reduction setting.) Under the assumption that $h$ is $L_h$-smooth and $g(x)=\tfrac{1}{m} \sum_{k=1}^{m} g_k(x)$ with $g_k(\cdot)$ being $L_{g_k}$-smooth, $\hat{L}_g = \max\limits_{k}L_{g_k}$
    and $\hat{L}_g \ge mL_h$, we propose an Accelerated Variance Reduction Method that obtains $\hat{x}$  such that $f(\hat{x}) - f(x^{\ast}) \leq \varepsilon$ with high probability in  $O\left(\sqrt{\tfrac{L_h}{\mu }} \ln\left( \tfrac{1}{\varepsilon} \right)\right)$ gradient evaluations for $h(\cdot)$ and $ O \left( \sqrt{\tfrac{m  \hat{L}_g}{\mu}} \ln\left( \tfrac{1}{\varepsilon} \right)  \right) $ stochastic gradient evaluations for $g(\cdot)$.
\end{enumerate}
In contrast to our results, in the full gradient setting, the standard approach \cite{monteiro2013accelerated,lin2015universal,lin2018catalyst} leads to the necessity to make $O\left(\sqrt{\tfrac{L_h+L_g}{\mu }} \ln\left( \tfrac{1}{\varepsilon} \right)\right)$ gradient evaluations \textit{both} for $h(\cdot)$ and $g(\cdot)$. Thus, in existing methods, each gradient is evaluated more times than it is suggested by our complexity bounds. Moreover, we can say that our bounds are in some sense optimal. Imagine that $g=0$ then by the lower complexity bound for smooth strongly convex optimization \cite{nemirovsky1983problem,nesterov2018lectures}, the number of gradient evaluations for $h$ can not be smaller than $O\left(\sqrt{\tfrac{L_h}{\mu }} \right)$, which coincides with our upper bound. A similar argument holds for $g$ in the case $h=0$.
\aai{We remark that the same full gradient setting for problem \eqref{prob_st} was considered in  \cite{lan2016accelerated}. In this particular case, our framework allows to obtain the same complexity results as in \cite{lan2016accelerated}, but by a different algorithm.}

\pd{
Further, existing coordinate descent methods do not allow to combine directly gradient oracle for $h$ and coordinate derivative oracle for $g$. If, in the random coordinate descent setting, one uses coordinate derivatives also for the function $h$, then the standard approach \cite{nesterov2017efficiency,gasnikov016accrand,allen2016even} leads to the necessity to make $O\left( \sum_{i=1}^n \sqrt{(\beta_i^h+\beta_i^g)/\mu} \ln\left( \tfrac{1}{\varepsilon} \right) \right)$ coordinate derivative evaluations \textit{both} for $h(\cdot)$ and $g(\cdot)$. Thus, in the existing methods, each coordinate derivative is evaluated more times than it is suggested by our complexity bounds. 
}

\pd{
Finally, existing variance reduced methods do not allow to combine directly the full gradient oracle for $h$ and stochastic gradient oracle for $g$. If, in the variance reduction setting, one reformulates  problem \eqref{prob_st} as 
\[
\min_{x \in \R^n}\tfrac{1}{m} \sum_{k=1}^{m} (h(x)+g_k(x)),
\]
then the standard approach \cite{allen2017katyusha,lan2018optimal,zhang2015stochastic,lin2014accelerated,lan2019unified} requires $ O \left( \sqrt{\tfrac{m  (L_h+\hat{L}_g)}{\mu}} \ln\left( \tfrac{1}{\varepsilon} \right)  \right) $ gradient evaluations for $h(\cdot)$ and $ O \left( \sqrt{\tfrac{m  (L_h+\hat{L}_g)}{\mu}} \ln\left( \tfrac{1}{\varepsilon} \right)  \right) $ stochastic gradient evaluations for $g(\cdot)$. Thus, in existing methods, each derivative is evaluated more times than it is suggested by our complexity bounds. 
Moreover, we can say that our bounds are in some sense optimal. 
For the number of gradient evaluations for $h$ the argument is the same as above for the case of full gradient oracles.
Imagine now that $h=0$ then by the lower complexity bound for smooth strongly convex finite-sum optimization \cite{agarwal2015lower}, the number of stochastic gradient oracle calls for $g$ can not be smaller than $ O \left( \sqrt{\tfrac{m  \hat{L}_g}{\mu}}  \ln\left( \tfrac{1}{\varepsilon} \right) \right) $, which coincides with our upper bound. 
}

\pd{We elaborate more on the comparison with the existing methods in Section \ref{exp}, where we consider two particular examples: Kernel Support Vector Machine (SVM) and regularized soft-max minimization motivated by Bayesian non-parametric approach to log-density estimation. 
}

\pd{Finally, as a part of our contribution and as a corollary of our general framework, we propose accelerated algorithms for the three settings described above without strong convexity assumption for $f$. In all the cases our framework is quite flexible and may use different accelerated methods as a building block to achieve oracle complexity separation. 
}

\section{Proposed Algorithmic Framework and Main Result}
\label{sec:1}
\aai{In this section we formally describe the proposed framework, state the main general theorem, which we apply in the next section in three particular settings, and give a sketch of the proof of the main result.}


\begin{algorithm}[htb]
\caption{\rev{Splitting Accelerated Envelope $\text{SAE}(x^0, L, N)$}}
\label{alg:MS}
\begin{algorithmic}[1]
\STATE {\bf Input:}  Starting point $x^0 = y^0 = z^0$; parameter $L \in (0, L_h]$; number of iterations $N$.
\FOR{$k = 0, 1, \ldots,N-1$}
        \STATE Compute 
           $a_{k+1} = \tfrac{1/L + \sqrt{1/L^2 + 4 A_k/L}}{2}$, $A_{k+1} = A_k + a_{k+1}$, $x^{k+1} = \tfrac{A_{k}}{A_{k+1}} y^{k} + \tfrac{a_{k+1}}{A_{k+1}} z^{k}$.
        \STATE \label{step:applyGMCO} \rev{Using the procedure GMCO outlined below, compute}  
	\begin{equation}
	    \label{aux_pr}
	    y^{k+1}  = \text{GMCO}(x^{k+1}, F_{L, x^{k+1}}(\cdot)).
	\end{equation}

    \STATE Compute  $z^{k+1}=z^{k}-a_{k+1} \nabla f\left(y^{k+1}\right)$.
    \ENDFOR
\STATE {\bf Output:} \aai{$\hat x := y^{N}$}.
\aai{\STATE {\bf procedure $\zeta^* = \text{GMCO}(\zeta^0, F_{L,\zeta^0}(\cdot))$} 
\STATE \label{step:GMCO_input}{\quad \bf Input:} starting point $\zeta^0 \in \mathbb{R}^n$, objective function $F_{L, \zeta^0}(\zeta) = f(\zeta) + \tfrac{L}{2} \|\zeta - \zeta^0\|^2_2$.
\STATE \quad {\bf for $k \geq 0$}
		\STATE  \label{step:GMCO_subproblem}  
		\begin{equation}\label{eq:GMCO_step}
		 \text{Set }	\varphi_{k}(\zeta)  :=  \left\langle \nabla h(\zeta^{k-1}), \zeta - \zeta^{k-1} \right\rangle + g(\zeta) + \tfrac{L}{2}\|\zeta-\zeta^0\|_2^2 + \tfrac{L_h}{2} \|\zeta - \zeta^{k-1}\|_2^2. 
		\end{equation}
		\STATE \label{eq:GMCO_step_1} \quad \quad  
		Run $\mathcal{M}_{inn}$ from Assumption \ref{ass_inn_method} using the number of $O_g$ calls $N_{\mathcal{M}}$ defined in ~\eqref{N_est}, and applied to the objective $\varphi_{k}(\zeta)$. Set  $\zeta^{k}:=\mathcal{M}_{inn}(\varphi_{k}(\zeta),\zeta^{k-1}, N_{\mathcal{M}})$.
		\STATE \label{step:checkMS} \quad \quad  \textbf{If} $\|\nabla F_{L, \zeta^0} (\zeta^k)\|_2 \leq \tfrac{L}{2} \|\zeta^k - \zeta^0\|_2 $, \textbf{then} break.
		\STATE \quad {\bf end for }
		\STATE \quad Output: $\zeta^* := \zeta^k$.
\STATE {\bf end procedure} }
\end{algorithmic}
\end{algorithm}

\rev{We start with a formal description of the proposed algorithmic framework in the convex case ($\mu=0$) and then describe a modification for the strongly convex case ($\mu>0$). Assumption \ref{ass_inn_method} governs the third loop, which is on the lowest level. Next, we describe algorithms that are used in the outer and middle level. We start with the first (outer loop), see Algorithm \ref{alg:MS}. The outer loop is based on 
the Monteiro--Svaiter Accelerated Proximal Method~\cite{monteiro2013accelerated} that
in each iteration requires to find an inexact proximal operator of the function $f$. This amounts to approximate minimization of the function
\begin{equation}
    \label{eq:F_L}
    F_{L,x}(y) :=  f(y) + \tfrac{L}{2}\|y - x\|^2_2,
\end{equation}
where the point $x$ is given and $L$ is a parameter that will be chosen later and is required to satisfy $0 < L \leq L_h$. We note also that directly applying the algorithm of \cite{monteiro2013accelerated} does not allow to separate the complexities as we desire. The reason is that in \cite{monteiro2013accelerated} the objective $f$ is considered as a whole and each iteration requires the same number of oracle calls both for $h$ and $g$.
}

\rev{In our setting, in each iteration of the outer loop, we need to approximately solve the following minimization problem
\[
\min\limits_{y} \left\{F_{L, \, x^{k+1}}(y)= h(y)+g(y)+
\frac{L}{2}\|y-x^{k+1}\|_2^2\right\}.
\]
This problem is $L$-strongly convex and has $L_h$-smooth part $h(y)$. Thus, it fits well the composite optimization framework \cite{nesterov2013gradient} with smooth part $h(y)$ and composite term $g(y) + \tfrac{L}{2}\|y - x^{k+1}\|^2_2$, and, to solve this problem, in step \ref{step:applyGMCO} of Algorithm \ref{alg:MS}, we apply  Gradient Method for Composite Optimization (GMCO)~\cite{nesterov2013gradient} listed as a separate procedure GMCO in  Algorithm~\ref{alg:MS}. This procedure constitutes the second/middle loop.
To make the connection with the third loop, we consider function $\varphi_k(\zeta)$ defined in \eqref{eq:GMCO_step} that needs to be minimized in each iteration of the procedure GMCO. 
This function has $L_g$-smooth part $g$ and strongly convex part $\tfrac{L}{2}\|\zeta-\hat{\zeta}^0\|_2^2 + \tfrac{L_h}{2} \|\zeta - \hat{\zeta}^{k-1}\|_2^2$, which can be considered as a composite term. Moreover, disregarding the constant terms, the structure of the function $\varphi_k(\zeta)$ in \eqref{eq:GMCO_step} is the same as that of function $\varphi(v)$ defined in \eqref{eq:innermost_problem} and, thus, is covered by Assumption \ref{ass_inn_method} with $\alpha = L + L_h$. Thus, to minimize $\varphi_k(\zeta)$ we can use as the inner method the method $\mathcal{M}_{inn}$ from Assumption~\ref{ass_inn_method}, where we take $\alpha = L + L_h$ and the number of iterations equal to 
\begin{equation}
\label{N_est}
    N_{\mathcal{M}} = 
  \begin{cases} 
  \tfrac{C_0\tau_{g}}{\sqrt{L + L_h}} \ln\tfrac{C_1 L_h R}{\delta \sqrt{\varepsilon L}}        & \text{if } \mu = 0, \\
   \tfrac{C_0\tau_{g}}{\sqrt{L + L_h}} \ln\tfrac{C_1 L_h }{\delta \sqrt{\mu L}}   & \text{if } \mu > 0,
  \end{cases}
\end{equation}
where $R \geq \|x^0 - x^{\ast}\|_2$ is an estimate for the distance to the solution of problem \eqref{prob_st}, $C_0, C_1>0$ are appropriately chosen numerical constants, and $\delta \in (0,1)$ is a desired probability confidence level.
}


\rev{In the strongly convex case ($\mu > 0$), the original algorithm \cite{monteiro2013accelerated} is not guaranteed to have   linear convergence. Thus, we provide the following extension of Algorithm \ref{alg:MS} using the restart technique. The result is listed below as Algorithm \ref{Alg:rest_str_MS}. }


\begin{algorithm}[ht]
\caption{Restarted Splitting Accelerated Envelope R-SAE}
\label{Alg:rest_str_MS}
\begin{algorithmic}[1]
\STATE {\bf Input:}  Starting point $\eta_0$; strong convexity constant $\mu > 0$; parameter $L > 0$.
\FOR{$t \geq 0$}
        \STATE Set $  N_0 = \sqrt{\tfrac{8L}{\mu}}$, find $\eta_{t}$ = $\text{SAE}(\eta_{t-1}, L, N_0)$.
    \ENDFOR
\STATE {\bf Output:} \aai{$\hat x := \eta_{t}$}
\end{algorithmic}
\end{algorithm}

\rev{Combining all the proposed loops together, we obtain that the complexity of our algorithmic framework is given by the following main theorem, which we prove in Section~\ref{App:proof}.
}
\begin{theorem}
\label{main_theo}
Let Assumptions 1 -- 4 hold and $L_g \geq L_h$. 
Then

\hspace{-1em}a) \rev{if $\mu = 0$, SAE (Algorithm \ref{alg:MS})} outputs a point $\hat{x}$  such that $f(\hat{x}) - f(x^{\ast}) \leq \varepsilon$ with probability at least $1 - \delta$ in 

\hspace{-2em}$O\left(\sqrt{\tfrac{LR^2}{\varepsilon}} \cdot \left(1 + \tfrac{L_h}{L} \right) \right)$ calls to the oracle $O_h$ for $h$ and 

\hspace{-2em}$O \left(\sqrt{\tfrac{LR^2}{\varepsilon}} \cdot \left(\kappa_g + \tfrac{L_h}{L}\cdot \left(\tfrac{\tau_{g}}{\sqrt{L + L_h}}\ln\tfrac{R}{\delta \sqrt{\varepsilon}} \right) \right) \right)$ calls to the basic oracle $O_g$ for $g$;

\hspace{-1em}b) \rev{if $\mu > 0$, R-SAE (Algorithm \ref{Alg:rest_str_MS})} outputs a point $\hat{x}$  such that $f(\hat{x}) - f(x^{\ast}) \leq \varepsilon$ with probability at least $1 - \delta$ in 

\hspace{-2em}$O\left(\sqrt{\tfrac{L}{\mu }} \ln\left(\tfrac{  \mu R^2}{\varepsilon} \right) \cdot \left(1 + \tfrac{L_h}{L} \right) \right)$ calls to the oracle $O_h$ for $h$ and
  
\hspace{-2em}$ O \left(\sqrt{\tfrac{L}{\mu}}  \ln\left(\tfrac{  \mu R^2}{\varepsilon}\right) \cdot \left(\kappa_g + \tfrac{L_h}{L}\cdot \tfrac{\tau_{g}}{\sqrt{L + L_h}}\ln\tfrac{1}{\delta }   \right) \right) $
  calls to the basic oracle $O_g$ for $g$.
\end{theorem}
    
  
  
  
As we see, we obtain separately the complexity estimates for the number of oracle calls for the function $h$ and the number of basic oracle calls for the function $g$. Moreover, we obtain these estimates in two settings: convex and strongly convex. 

Let us give a proof sketch in a simple setting with $\mu >0$ and $O_g$ producing the full gradient $\nabla g(\cdot)$.
In the outer loop we solve problem \eqref{prob_st} by accelerated proximal method \cite{ivanova2019adaptive}. This method takes as an input parameter $L$ satisfying $\mu \le L \le L_h$ and the total number of its iterations to find an approximate solution for problem \eqref{prob_st} is $\tilde{O}(\sqrt{L/\mu})$.\footnote{ Here and below, for simplicity, we hide numerical constant and polylogarithmic factors using non-asymptotic $\tilde{O}$-notation. More precisely, $\psi_1(\varepsilon,\delta) = \tilde{O}(\psi_2(\varepsilon,\delta))$ if there exist constants $C,a,b>0$ such that, for all $\varepsilon>0$, $\delta \in (0,1)$, $\psi_1(\varepsilon,\delta) \leq C\psi_2(\varepsilon,\delta)\ln^a\frac{1}{\varepsilon}\ln^b\frac{1}{\delta}$.}
Moreover, each iteration of this algorithm requires to find, possibly inexactly, a proximal operator of $f$, i.e. to solve the problem
\begin{equation}\label{composite}
\min_{x\in\R^n} h(x) + g(x) + \frac{L}{2}\|x-x^{k+1}\|^2_2,
\end{equation}
where $x^{k+1}$ is the current iterate. The latter problem is solved in the next loop by non-accelerated composite gradient method \cite{nesterov2013gradient} with $g(x) + \frac{L}{2}\|x-x^{k+1}\|^2_2$ being the composite term. This is an $L$-strongly convex problem with $L_h$-smooth part $h(x)$, and each of the $\tilde{O}(L_h/L)$ iterations of this method 
requires to solve a problem with the structure $\varphi(v) = \left\langle \beta, v \right\rangle + \tfrac{L + L_h}{2}\|v\|_2^2 + g(v)$, which, for $\alpha=L + L_h$ is exactly the structure given in Assumption \ref{ass_inn_method}. Moreover, this problem is $(L+L_h)$-strongly convex and has $L_g$-smooth part $g$. Thus, we can solve this problem by accelerated gradient method for composite optimization \cite{nesterov2013gradient} in $\tilde{O}(\sqrt{L_g/(L+ L_h\att{)}})$ evaluations of $\nabla g(\cdot)$ (we have $\tau_g=\sqrt{L_g}$). This is the last loop and we can estimate the total number of calls for the gradient of each function $h$ and $g$. The gradient of $h$ is used in the second loop, and is called in total $\tilde{O}(\sqrt{L/\mu})\cdot \tilde{O}(L_h/L)$ times. The gradient of $g$ is used in the third loop, and is called in total
\begin{equation}\label{complexity}
    \tilde{O}(\sqrt{L/\mu})\cdot \left[\tilde{O}(L_h/L)\cdot \tilde{O}\left(\sqrt{L_g/(L+L_h)}\right)\right]
\end{equation}
times. Minimizing this expression in $ L \in\left[\mu, L_h\right]$, 
we obtain that we should choose $L = L_h$.
Thus, we can solve problem \eqref{prob_st} via
\begin{center}
    $\tilde{O}\left(\sqrt{L_h/\mu}\right)$ evaluations of $\nabla h(x)$  \;\; and \;\; $\tilde{O}\left(\sqrt{L_g/\mu}\right)$ evaluations of $\nabla g(x)$.
\end{center}
Note that we can solve the optimization problem in the innermost loop by non-composite accelerated gradient method \cite{nesterov2018lectures}. Then the number of evaluations of $\nabla g(\cdot)$ in this loop is $\tilde{O}(\sqrt{(L_g+L+ L_h)/(L+ L_h\att{)}})$, but under the assumption that $L_g\geq L_h$ we still obtain the same complexity bounds for the number of oracle calls for $h$ and $g$.
For simplicity we have simplified here part of the derivations, and the full proofs, including the case of randomized oracles and non-strongly convex optimization, are given below. One of our main technical contributions is the careful analysis of how the number of iterations in each loop affects the (probabilistic) inexactness of the solution to the auxiliary optimization problem in the next level and eventually the complexity of the whole method.

\vspace{-1em}

\section{Applications}\label{appl}
\pd{In this section, we consider three particular application examples for our general framework and general Theorem \ref{main_theo}. By considering three different types of basic oracle $O_g$ for $g$: full gradient, coordinate derivative, stochastic gradient in the finite-sum setting, and considering three corresponding different methods as $\mathcal{M}_{inn}$, we obtain accelerated gradient method, accelerated random coordinate descent and accelerated stochastic variance reduced method with oracle complexity separation.}
\vspace{-1em}

\subsection{Accelerated Gradient Method}
Assume that for the problem \eqref{prob_st} $L_g \ge L_h$ and the basic oracle $O_g$ returns the full gradient of the function $g$, meaning that $\kappa_g = 1$. In this case, we have two options for solving the problem \eqref{eq:GMCO_step} in the third loop and ensuring Assumption \ref{ass_inn_method}, for which we have $\alpha=L+L_h$.
Firstly, considering $\left\langle \beta, v \right\rangle + \tfrac{\alpha}{2}||v||_2^2$ in $\varphi(v)$ in Assumption \ref{ass_inn_method} as a composite term, we can apply Accelerated Gradient Method for Composite Optimization from~\cite{nesterov2013gradient} as $\mathcal{M}_{inn}$ in the third loop, which means that Assumption~\ref{ass_inn_method} holds with $\tau_g = \sqrt{L_{g}}$ since the smooth part $g$ of the objective is $L_g$-smooth. Secondly, we can apply Accelerated Gradient Method \cite{nesterov2018lectures} as $\mathcal{M}_{inn}$ to the whole objective in the third loop, which means that Assumption~\ref{ass_inn_method} holds with $\tau_g = \sqrt{L_{g}+L+L_h}$ since in this case the objective is $(L_{g}+L+L_h)$-smooth.
In both cases we can choose $L=L_h$ and obtain  $\tau_g =O(\sqrt{L_{g}})$ under our assumption that $L_g \ge L_h$.  
This leads to the following corollary of Theorem~\ref{main_theo}.

\begin{corollary}
    Under assumption $L_g \ge L_h$, using the Accelerated Gradient Method (for Composite Optimization) as $\mathcal{M}_{inn}$, \rev{respectively SAE and R-SAE} obtain $\hat{x}$  such that $f(\hat{x}) - f(x^{\ast}) \leq \varepsilon$ in the following number of oracle calls:
    
  \hspace{-1em}a) $O\left(\sqrt{\tfrac{L_\ag{h} R^2}{\varepsilon}} \right)$ calls for $h(\cdot)$, $O \left(\sqrt{\tfrac{L_g R^2}{\varepsilon}} \right)$ calls for $g(\cdot)$, if $\mu = 0$, and  
  
  \hspace{-1em}b)\aai{ $O\left(\sqrt{\tfrac{L_\ag{h}}{\mu }} \ln\tfrac{ 1}{\varepsilon} \right)$  calls for $h(\cdot)$ and $ O \left(\sqrt{\tfrac{L_g}{\mu}} \ln\tfrac{ 1}{\varepsilon}  \right) $  calls for $g(\cdot)$, if $\mu > 0$.}
\end{corollary}
\aai{We remark that the same full-gradient oracles setting for problem \eqref{prob_st} was considered in  \cite{lan2016accelerated}. In this particular case, our framework allows to obtain the same complexity results as in \cite{lan2016accelerated}.}
\vspace{-1em}

\subsection{Accelerated Random Coordinate Descent Method}
\label{S:ARCD}
Assume that in the problem \eqref{prob_st} the part $g(\cdot)$ has block smoothness, i.e. there exist coordinate Lipschitz constants $\beta_1, \dots, \beta_n$ such that for any $ x \in \mathbb{R}^n, u \in \mathbb{R}$
$
\left|\nabla_{i} g\left(x+u e_{i}\right)-\nabla_{i} g(x)\right| \leq  \beta_{i}|u|, \quad  i =1, \ldots, n,
$
where $\nabla_{i} g(x) = \partial g(x) / \partial x_i$ and $e_i$ is the $i$-th coordinate vector.
For twice differentiable function $g(\cdot)$ this condition is equivalent to the condition $(\nabla^{2} g(x))_{i, i} \leq \beta_{i}.$
As the basic oracle $O_g$ we take an oracle \ai{computing the partial derivative $\nabla_i g(\cdot)$ for a given $i$}. Thus, we need $\kappa_g = n$ calls to $O_g$ to compute the full gradient $\nabla g(\cdot)$.
In this case, we apply accelerated random coordinate descent as $\mathcal{M}_{inn}$ in the third loop to solve problem \eqref{eq:GMCO_step} and ensure Assumption \ref{ass_inn_method}, for which we have $\alpha=L+L_h$.
As before, we have two options. Firstly, considering $\left\langle \beta, v \right\rangle + \tfrac{\alpha}{2}||v||_2^2$ in $\varphi(v)$ in Assumption \ref{ass_inn_method} as a composite term, we can apply Accelerated Proximal Random Coordinate Method from~\cite{fercoq2015accelerated,allen2016even} as the inner method $\mathcal{M}_{inn}$. For this method Assumption~\ref{ass_inn_method} holds with $\tau_g =  n \sqrt{\overline{L}_{g}},$ where $\sqrt{\overline{L}_{g}} = \tfrac{1}{n}  \sum_{i=1}^n \sqrt{\beta_i}$. Secondly, we can apply Accelerated Random Coordinate Method~\cite{nesterov2012efficiency,nesterov2017efficiency,gasnikov016accrand} as the inner method $\mathcal{M}_{inn}$ for the whole objective in \eqref{eq:GMCO_step}. Then Assumption~\ref{ass_inn_method} holds with $\tau_g =  n \sqrt{\overline{L}},$ where $\sqrt{\overline{L}} = \tfrac{1}{n}  \sum_{i=1}^n \sqrt{\beta_i+L+L_h}$. 
Let us assume that $\overline{L}_{g} \ge L_h$. Then in both cases we can choose $L=L_h$ and obtain  $\tau_g =O(n \sqrt{\overline{L}_{g}})$.  
Summarizing, we obtain the following corollary of Theorem~\ref{main_theo}.

\begin{corollary}\label{coor}
    Under assumption $\overline{L}_{g} \ge L_h$, using the Accelerated (Proximal) Random Coordinate Method as $\mathcal{M}_{inn}$ \rev{respectively SAE and R-SAE} obtain $\hat{x}$  such that $f(\hat{x}) - f(x^{\ast}) \leq \varepsilon$  with probability at least $1 - \delta$  in the following number of oracle calls:
    
 \hspace{-2em}a) $O\left(\sqrt{\tfrac{L_h R^2}{\varepsilon}} \right)$  calls for $h(\cdot)$, $O \left(n  \sqrt{\tfrac{\overline{L}_{g} R^2}{\varepsilon}}\ln\frac{1}{\delta} \right)$  calls for $g(\cdot)$, if $\mu = 0$, and   
  
 \hspace{-2em}b) \aai{$O\left(\sqrt{\tfrac{L_h}{\mu }}\ln \tfrac{1}{\varepsilon}   \right)$  calls for $h(\cdot)$ and $ O \left(n  \sqrt{\tfrac{\overline{L}_{g}}{\mu}}   \ln \tfrac{1}{\varepsilon} \cdot \ln\frac{1}{\delta} \right) $  calls for $g(\cdot)$, if $\mu > 0$.}

\end{corollary}

\ag{Note, that if $\mathcal{M}_{inn}$ is a directional search or a derivative-free method such as in \cite{dvurechensky2017randomized}, then the main conclusions of Corollary~\ref{coor} remain valid after replacing $\overline{L}_{g}$ with $L_{g}$.}

\subsection{Accelerated Stochastic Variance Reduced Method}
\label{S:ASVRG}
Consider the following minimization problem
\begin{equation}
\label{eq_h_sumg}
\min_{x\in \R^n} \left\{ f(x)= h(x) + \tfrac{1}{m} \sum_{k=1}^{m} g_k(x) \right\}
\end{equation}
that is \eqref{prob_st} with $g(x)=\tfrac{1}{m} \sum_{k=1}^{m} g_k(x)$.
We assume that each component $g_k(\cdot)$ is $L_{g_k}$-smooth. As the basic oracle $O_g$ we take an oracle that, given $k$, computes  $\nabla g_k(\cdot)$. Thus, we need $\kappa_g = m$ basic oracle $O_g$ calls to compute the full gradient $\nabla g(\cdot).$ 
In this case, we would like to apply accelerated stochastic variance reduced algorithms as $\mathcal{M}_{inn}$ to solve problem \eqref{eq:GMCO_step} in the third loop and ensure Assumption \ref{ass_inn_method}, for which we have $\alpha=L+L_h$.
As before, we have two options. Firstly, considering $\left\langle \beta, v \right\rangle + \tfrac{\alpha}{2}\|v\|_2^2$ in $\varphi(v)$ in Assumption \ref{ass_inn_method} as a composite term, we can apply composite versions of Accelerated Stochastic Variance Reduced Algorithms, e.g. Katyusha~\cite{allen2017katyusha} or Varyag \cite{lan2019unified} as the inner method $\mathcal{M}_{inn}$. 
For these methods the number of oracle calls to solve problem \eqref{eq:GMCO_step} is $\tilde{O}\left(m + \sqrt{\tfrac{m \hat{L}_g}{L}} \right)$, where $\hat{L}_g = \max\limits_{k}L_{g_k}$, and Assumption~\ref{ass_inn_method} holds with $\tau_g =  \sqrt{m \hat{L}_g}$.
Secondly, we can apply Accelerated Stochastic Variance Reduced Algorithms \cite{allen2017katyusha,lan2018optimal,lan2019unified} as the inner method $\mathcal{M}_{inn}$ for the whole objective in \eqref{eq:GMCO_step}. Then Assumption~\ref{ass_inn_method} holds with $\tau_g =  \sqrt{m (\max\limits_{k}L_{g_k}+L+L_h)}$.
Let us assume that $mL_{h}  \leq \hat{L}_g$. Then in both cases we can choose $L=L_h$ and obtain $\tau_g =  O\left(\sqrt{m \hat{L}_g}\right)$.
As a result, we obtain the following corollary of Theorem~\ref{main_theo}.

\begin{corollary}
    Under assumption that $ mL_{h} \leq \hat{L}_g$, using Accelerated (Composite) Stochastic Variance Reduced algorithm as $\mathcal{M}_{inn}$ \rev{respectively SAE and R-SAE} obtain $\hat{x}$  such that $f(\hat{x}) - f(x^{\ast}) \leq \varepsilon$  with probability at least $1 - \delta$ in the following number of oracle calls:
    
  \hspace{-2em}a)  $O\left(\sqrt{\tfrac{L_\ag{h} R^2}{\varepsilon}} \right)$  calls for $h(\cdot)$, $O \left( \sqrt{\tfrac{m \hat{L}_g R^2}{\varepsilon}} \right)$  calls for $g(\cdot)$, if $\mu = 0$, and  
  
  \hspace{-2em}b) \aai{$O\left(\sqrt{\tfrac{L_h}{\mu }}\ln \tfrac{1}{\varepsilon}    \right)$  calls for $h(\cdot)$ and $ O \left( \sqrt{\tfrac{m  \hat{L}_g}{\mu}}   \ln \tfrac{1}{\varepsilon}\cdot \ln \tfrac{1}{\delta}     \right) $  calls for $g(\cdot)$, if $\mu > 0$.}

\end{corollary}

\section{Experiments}\label{exp}
In this section, we investigate the practical performance of our algorithmic framework on two particular machine learning problems: Kernel Support Vector Machine (SVM) problem and log-density estimation using Bayesian approach. We also give some theoretical explanations for the observed results.

\subsection{Kernel Support Vector Machine}\label{SVM}
We start with one of the basic machine learning problems: a binary classification problem using the Kernel trick  \cite{NatureStLearning,SVM_book}. Given a set of feature vectors $a_k$, $k=1,...,m$ 
and known classes $b_k \in  \left\lbrace -1, \, +1\right\rbrace$, $k=1,...,m$, the goal is to find linear classifier by solving the following Kernel Support Vector Machine (SVM) problem
\begin{equation}
 \label{exp_eq_svm_main}
 \min_{x_0, \, x \in \R^n} \frac{1}{m} \sum_{k=1}^m(1-\langle b_k,x_0+Kx\rangle)_+ + \frac{\lambda}{2}\langle x,Kx\rangle,
\end{equation}
where 
$[K]_{ij}=K(a_i,a_j)$, $i,j=1,\ldots,n$ is 
some chosen positive definite kernel matrix
and $\lambda >0$ is the regularization parameter. 
 

Since our general framework is constructed for smooth minimization problems, we apply Nesterov's smoothing technique~\cite{nesterov2005smooth},
as described in \cite{Zhang_regularizedrisk}, and change the hinge loss to its smooth softmax approximation. This gives us a minimization problem of the form~\eqref{eq_h_sumg} with a quadratic term $h$ \rev{and smoothing parameter $\mu>0$}:
\begin{equation}
 \label{exp_eq_svm_main_smooth}
 \rev{\min_{x_0, \, x} \frac{1}{m} \sum_{k=1}^m \mu \ln \left(1+\exp \left(\frac{1}{\mu}(1-\langle b_k,x_0+Kx\rangle)\right) \right)+ \frac{\lambda}{2}\langle x,Kx\rangle.}
\end{equation}
As we consider a large-dimensional problem, the proximal operator for the quadratic term is too expensive since it requires inversion of the kernel matrix~$K$. Hence, as  outlined in Section \ref{S:ASVRG}, we use full gradient oracle for~$h(x)=\frac{\lambda}{2}\langle x,Kx\rangle$ and stochastic gradient oracle for the sum of smoothed hinge losses. 
In the experiments we use variance-reduction-based method Katyusha~\cite{allen2017katyusha} as the inner method $\mathcal{M}_{inn}$ in Algorithm~\ref{alg:MS}.

For the experiments, we chose publicly available
"Gas sensors for home activity monitoring" dataset\footnote{\url{http://archive.ics.uci.edu/ml/datasets/gas+sensors+for+home+activity+monitoring}} containing in our case 6200 instances of 11-dimensional feature vectors.
This dataset has recordings of a gas sensor array composed of 8 gas sensors, and temperature and humidity sensors. This sensor array was exposed to background home activity while subject to two different stimuli. The aim is to learn how to discriminate among them.
As a kernel function
we use the radial kernel $K(a, a') = \exp(-\gamma \|a - a'\|_2^2)$ with different values of~$\gamma$. \rev{Importantly, in this case, kernel $K$ turned out to be very ill-conditioned with the smallest eigenvalue as small as $10^{-16}$. Thus, essentially, the objective in \eqref{exp_eq_svm_main_smooth} is not strongly convex.}

\vspace{-1em}
\begin{figure}[ht]

    \begin{center}
        \includegraphics[width=0.75\textwidth]{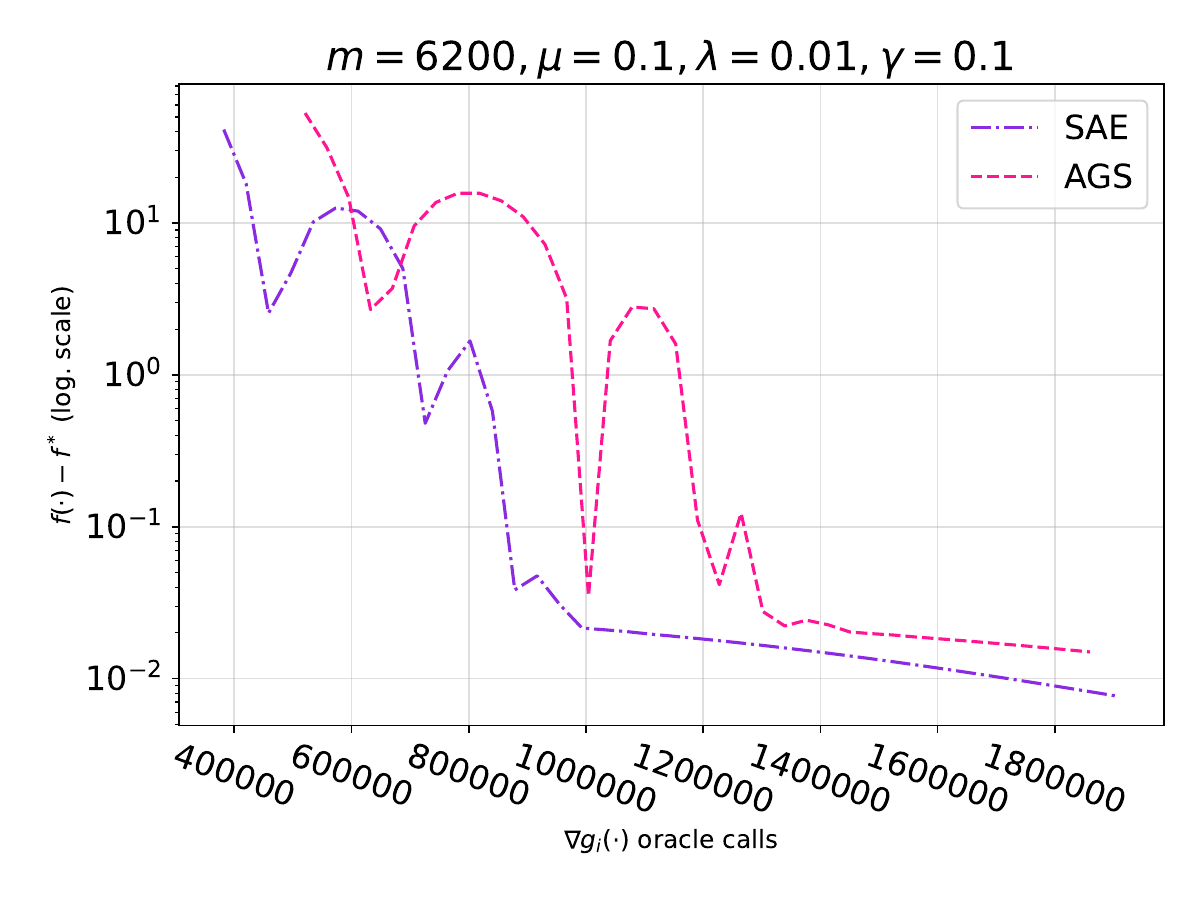}
        \caption{Performance of SAE (Algorithm~\ref{alg:MS}) and AGS \cite{lan2016accelerated} in terms of the number of $\nabla g_i(\cdot)$ oracle calls.}
        \label{svm_1}
    \end{center}
\end{figure}

\begin{figure}[ht]
    \begin{subfigure}[b]{0.5\textwidth}
        \includegraphics[width=\textwidth]{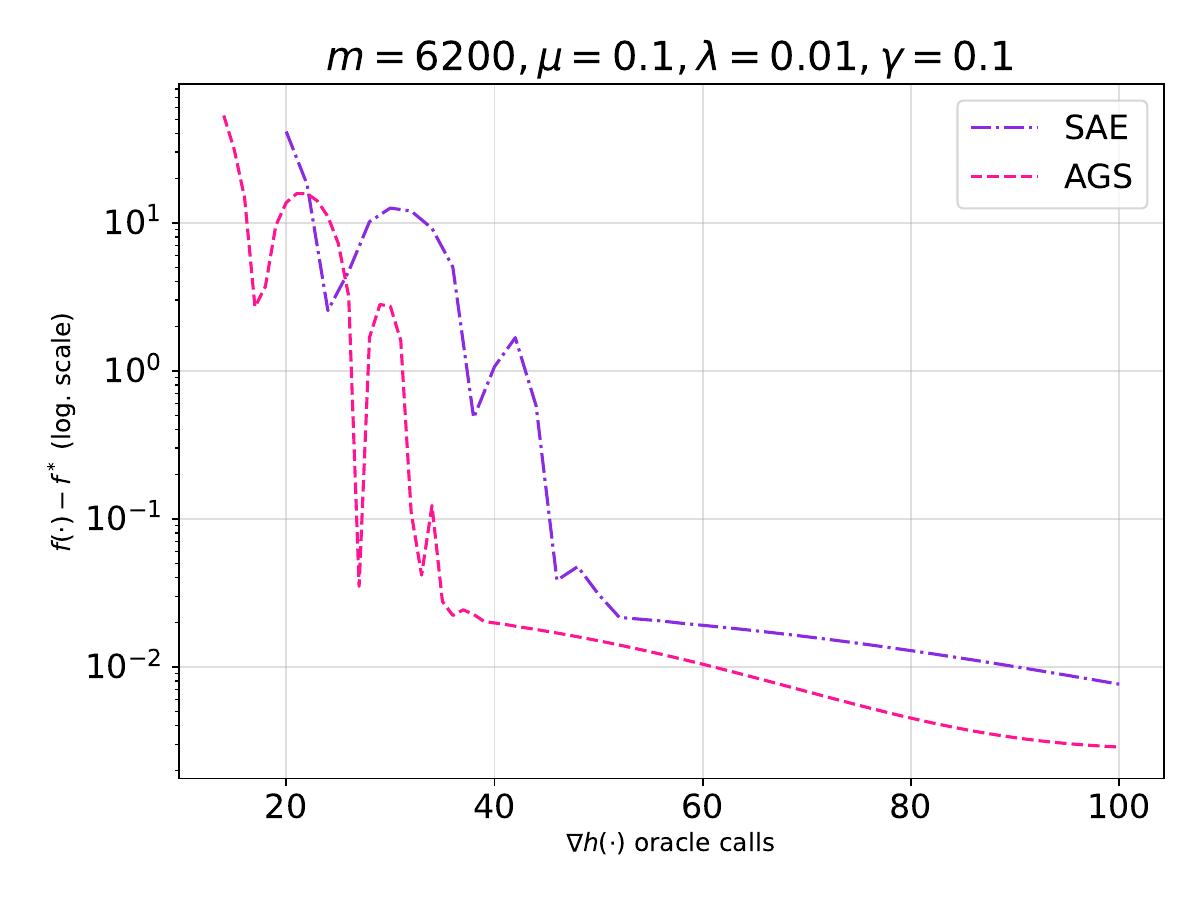}
        \caption{Performance of SAE (Algorithm~\ref{alg:MS}) and AGS \cite{lan2016accelerated} in terms of the number of $\nabla h(\cdot)$ oracle calls.}
    \end{subfigure}
    \hspace{1em}
    \begin{subfigure}[b]{0.5\textwidth}
        \includegraphics[width=\textwidth]{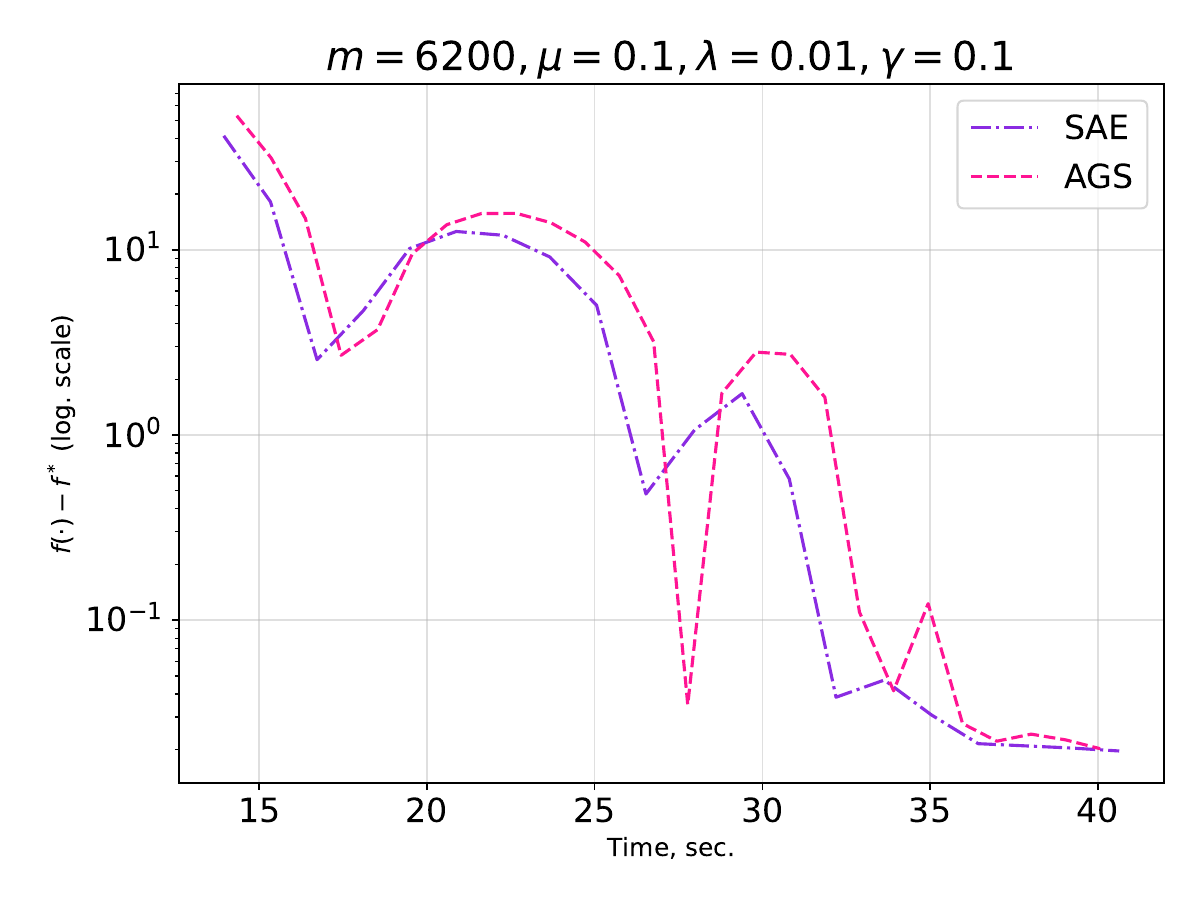}
        \caption{Performance of SAE (Algorithm~\ref{alg:MS}) and AGS \cite{lan2016accelerated} in terms of wall-clock  time. \newline}
    \end{subfigure}
    \caption{Comparison of SAE (Algorithm~\ref{alg:MS}) and AGS~\cite{lan2016accelerated}}
    \label{svm_2}
    
\vskip -0.2in
\end{figure}

Figures~\ref{svm_1} and \ref{svm_2} summarize the results of the experiments\footnote{Source code of these experiments is available at: \url{https://github.com/dmivilensky/Sliding-for-Kernel-SVM}}.
We compared \rev{the performance of the SAE (Algorithm~\ref{alg:MS})} with the \rev{performance of
Accelerated Gradient Sliding (AGS) method, proposed in \cite{lan2016accelerated}
}.
The results of the experiments show that
Algorithm~\ref{alg:MS} allows one \rev{to call the stochastic gradient oracle $g_i(\cdot)$ less often, which affects the working time of the method and allows SAE to work faster than AGS.}


The experiments reported in this Section
were performed on a personal computer with an 6 Intel Core i5 3 GHz CPU and 8 GB of RAM, under 64-bit macOS operating system, with the Jupyter Notebook browser environment. The version of Python programming language is 3.8.2, IPython 7.13.0.

Next we compare  theoretical complexity bounds for our approach, AGS \cite{lan2016accelerated} and Fast Gradient Method (FGM)  \cite{nesterov2018lectures} \rev{, by rewriting \eqref{exp_eq_svm_main_smooth} as follows}
\rev{
    $$\min_{x\in\R^n} \frac{1}{2}\langle x, Cx \rangle + \frac{1}{m}\sum_{k=1}^m f_k\left(\langle A_k,x\rangle \right),$$
}
where, to simplify the notation we denote $C=\lambda K$ and $A_k=K^Tb_k$, $k=1,\ldots,m$. 
\rev{In the experiments we observed that the smallest eigenvalue of the kernel $K$ is very small (and, hence, the same holds for $C$). Thus, we consider the above problem under the assumption of convexity, rather than strong convexity.}
Since we use the smoothing technique, we have
\rev{$|f_k''(y)| = O(1/\varepsilon)$}, $k=1,\ldots,m$, where $\varepsilon>0$ is the desired accuracy of the solution.
We assume that matrix $A = [A_1,...,A_m]^{\top}$ has $ms$ nonzero elements and that $\max_{k=1,...,m} \|A_k\|_2^2 = O(s)$, where $1\ll s \le n$ and $C$ is positive semidefinite  matrix 
with $\lambda_{\max}(C)\le 1/(\varepsilon m)$. We also assume that for the starting point $x_0$ and the closest solution $x^*$ it holds that $\|x_0-x^*\|_2\leq R$. In this setting, 
FGM \cite{nesterov2018lectures} requires $O\left(\sqrt{\frac{\left(s/\varepsilon + \lambda_{\max}(C)\right)R^2}{\varepsilon}}\right)$
iterations with the complexity of each iteration being
$O\left(ms + n^2\right)$ arithmetic operations, where the first term is obtained if one uses a special representation of matrix $A$ as an adjacency list.
Using our approach it is sufficient to make 
$\tilde{O}\left(\sqrt{\frac{ \lambda_{\max}(C)R^2}{\varepsilon}}\right)$ calls to the full gradient oracle $\rev{\nabla}h$ with the complexity of each call $O(n^2)$ and $\tilde{O}\left(\sqrt{\frac{\left(ms/\varepsilon \right)R^2}{\varepsilon}}\right)$ calls to the stochastic gradient oracle \rev{using the variance reduced algorithm \cite{allen2017katyusha}} with the complexity of each call being $O(s)$ arithmetic operations. \rev{Similarly, AGS algorithm \cite{lan2016accelerated}  requires to evaluate $\nabla h$ oracle  $O\left(\sqrt{\frac{ \lambda_{\max}(C)R^2}{\varepsilon}}\right)$ times with the complexity of each call $O(n^2)$,  and the full gradient $\nabla g$ (non-stochastic) oracle $O\left(\sqrt{\frac{\left(s/\varepsilon \right)R^2}{\varepsilon}}\right)$ times, while the complexity of the latter is $O(ms)$ arithmetic operations}.
These results are summarized in Table 1 below. The main observation is that when $s\gg 1$ and $\lambda_{\max}(C)\le 1/(\varepsilon m)\ll s/\varepsilon$ our approach leads to \rev{theoretical complexity that is better than for FGM}. \rev{Moreover, the complexity of the proposed approach with respect to the number $m$ is asymptotically better than for AGS ($O(m)$ in AGS vs $\tilde{O}(\sqrt{m})$ in SAE)}.
\begin{table}[ht]
\label{table_0}
\begin{center}
\begin{tabular}{|c|c|c|} \hline
    Algorithm  & Complexity &  Reference \\ \hline
    FGM        & $O\left(\frac{R}{\varepsilon}\sqrt{s}\left(ms + n^2\right)\right)$ & \cite{nesterov2018lectures}\\ \hline
    \rev{AGS}  & $O\left(\frac{R}{\varepsilon}\sqrt{s}\cdot ms\right) + O\left(\sqrt{\frac{\lambda_{\max}(C)R^2}{\varepsilon}}\cdot n^2\right)$ & \cite{lan2016accelerated}\\ \hline
    \rev{SAE}  & $\tilde{O}\left(\frac{R}{\varepsilon}\sqrt{sm}\cdot s\right) + \tilde{O}\left(\sqrt{\frac{\lambda_{\max}(C)R^2}{\varepsilon}}\cdot n^2\right)$ & this paper\\ \hline
\end{tabular}
\caption{Comparison of theoretical complexity bounds for SAE (Algorithm~\ref{alg:MS}), AGS and FGM.}
\end{center}
\end{table}
%
%
\subsection{Log-Density Estimation with Gaussian Prior}\label{subs_log}
In this subsection we consider log-density estimation problem \cite{spokoiny2019accuracy} for some unknown probability measure $\mathbb{P}$ using Gaussian prior using Bayesian approach. We assume that we have $m$ random observations $\tilde{z}_1,\dots, \tilde{z}_m \in \Z$ sampled from the measure $\mathbb{P}$. Without loss of generality, we assume that $\Z$ has finite support $\left\lbrace   z_k\right\rbrace_{k=1}^p$ of size $p$, whence
$\sum_{k=1}^p f(z_k) =1$, where $f$ is the unknown density function for the measure $\mathbb{P}$.
We parameterize the log-density by the linear model, i.e. we assume that $\ln f(z) = \sum_{i=1}^n x_i^{\ast}a_i(z) - c(x^{\ast})$,
where $a_1(z), a_2(z), \dots, a_n(z)$ are given basis functions and $x^{\ast} \in \R^n$ is an unknown vector, corresponding to the true density. The normalization constant $c(x^{\ast})$ is defined by 
\begin{equation}
    c(x) = \ln\left( \sum_{k=1}^p\exp\left(\left\langle A_k, x\right\rangle\right) \right),
\end{equation}
where  $A_k = a(z_k) = (a_1(z_k), \dots, a_n(z_k))^{\top}$ is the $k$-th column of $A = [a_j(z_k)]_{j,k = 1}^{n,p}$.
It is known from the Fisher theorem, see, e.g. \cite{spokoiny2019accuracy}, that the true parameter $x^{\ast}$ can be estimated via Maximum Likelihood Estimation (MLE) 
\[
\tilde x = \arg\max_{x\in \R^n} \left\lbrace \sum_{k=1}^m\left\langle  a(\tilde{z}_k), x \right\rangle - mc(x) \right\rbrace.
\]
Moreover, if we introduce a Gaussian prior $\mathcal{N}(0,G^2)$ for $x^{\ast}$, MLE translates to solving the following optimization problem
\begin{equation}\label{bayes}
\tilde x_G = \arg\max_{x\in \R^n} \left\lbrace \sum_{k=1}^m\left\langle  a(\tilde{z}_k), x \right\rangle - mc(x) -\frac{1}{2}\|Gx\|^2\right\rbrace.
\end{equation}
Bernstein--von Mises theorem claims \cite{spokoiny2019accuracy}, that $\tilde x_G$ is a good estimate of $x^{\ast}$ in the Bayesian approach.
For the numerical experiments we consider a particular case when matrix $A$ is sparse and all elements of $G^2$ are from the interval [1,2]. Modern Accelerated Random Coordinate  Descent (ARCD) algorithms \cite{fercoq2015accelerated} do not allow to take into account the sparsity of matrix $A$, so for the first two terms under $\argmax$ of r.h.s. of \eqref{bayes} it would be better to use standard Fast Gradient Method (FGM)  \cite{nesterov2018lectures}. The third term in \eqref{bayes} is, on the contrary, very friendly for applying  ARCD \cite{nesterov2017efficiency}. Under these assumptions we solve problem \eqref{bayes} with relatively small $m$ (or relatively large $G^2$) by our oracle complexity separation framework with $\mathcal{M}_{inn}$ being ARCD as described in Section \ref{S:ARCD}. 

Motivated by the above example, we solve the following problem:
\begin{equation}\label{eq:softmax_reg}
\min_{x \in \mathbb{R}^n} \left\{f(x) = \ln\left(\sum_{k=1}^p \exp \left(\left\langle A_k, x \right\rangle\right)\right) + \frac{1}{2} \|G x\|_2^2\right\}, 
\end{equation}
\rev{defining $h(x) = \ln\left(\sum_{k=1}^p \exp \left(\left\langle A_k, x \right\rangle\right)\right)$, $g(x)=\frac{1}{2} \|G x\|_2^2$.}
In our case, $n=500$, $p=6000$, $A$ is a sparse $p \times n$ matrix with sparsity coefficient $0.001$ and non-zero elements drawn randomly from $\mathcal{U}(-1, 1)$, and matrix $G^2$ is generated as $G^2 = \sum_{i=1}^n \lambda_i \tilde{e}_i^\top \tilde{e}_i$,
where $\sum_{i=1}^n \lambda_i = 1$ and $\left[\tilde{e}_i\right]_j \sim \mathcal{U}(1, 2)$ for every $i, j$. The gradient Lipschitz constant for the first term of $f$ is estimated as $L_h = \max_{i = 1,...,n} \|A^{\left\langle k\right\rangle}\|_2^2$,
where $A^{\left\langle k\right\rangle}$ is the $k$-th column of $A$. We take the parameter $L = 25 L_h$. Then  Lipschitz constants for the directional derivatives of the function $\varphi$ in \eqref{eq:innermost_problem} 
are estimated as $\beta_i = G^2_{ii} + L + L_h$. 

\begin{figure}[ht]
\begin{subfigure}[b]{0.5\textwidth}
\includegraphics[width=\columnwidth]{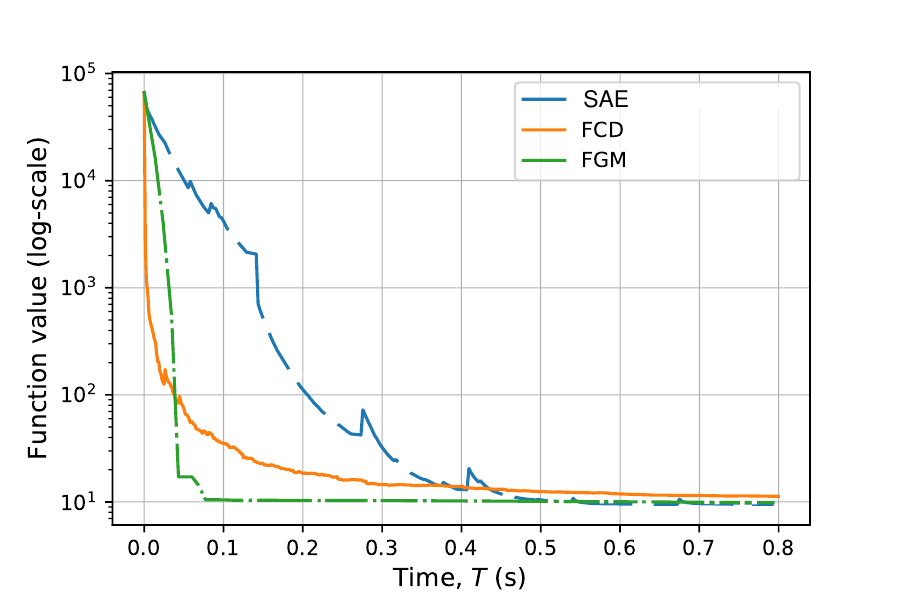}
\caption{Starting 0s}
\label{coord1}
\end{subfigure}
\begin{subfigure}[b]{0.5\textwidth}
\includegraphics[width=\columnwidth]{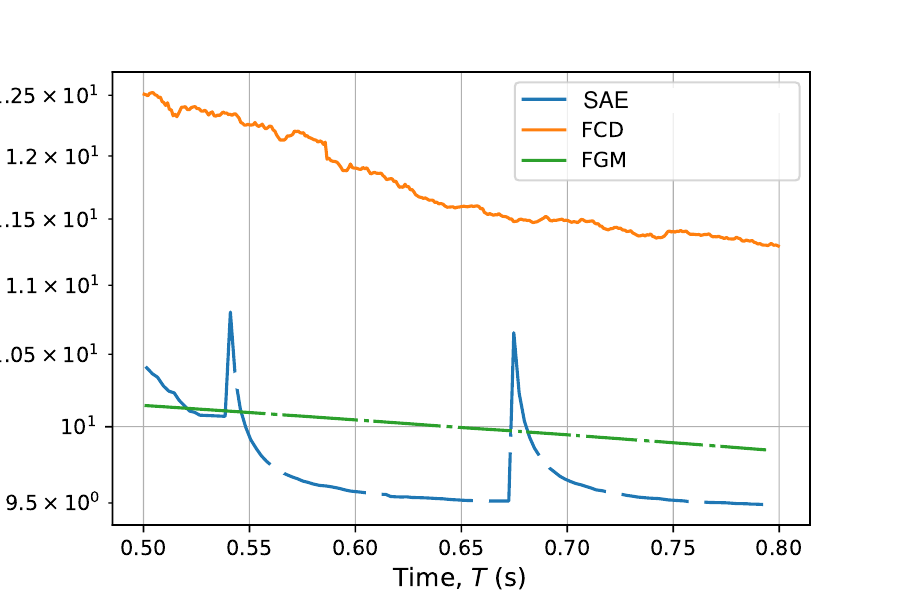}
\caption{Starting 0.5s}
\label{coord2}
\end{subfigure}
\caption{Objective value $f(\rev{y}^k)$ vs working time}
\vskip -0.2in
\end{figure}

\begin{figure}[ht]
\vskip 0.2in
\begin{center}
\centerline{\includegraphics[width=\columnwidth]{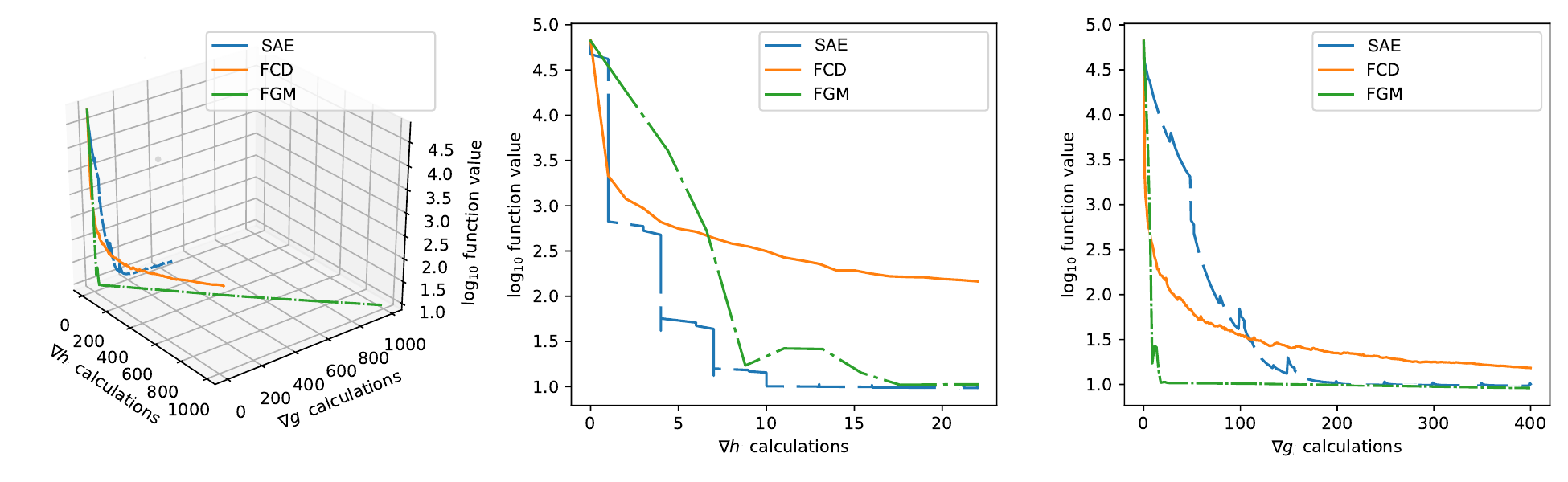}}
\caption{Objective value $f(\rev{y}^k)$ vs number of $\nabla h(\cdot)$ and $\nabla g(\cdot)$ calls}
\label{grad_calc}
\end{center}
\vskip -0.2in
\end{figure}

Following Section \ref{S:ARCD}, as $\mathcal{M}_{inn}$  we use  restarted every $300$ iterations ARCD with $\beta=1/2$ \cite{nesterov2017efficiency}. We compare our approach with ARCD (denoted as FCD in the figures) and FGM applied to the whole objective without complexity separation. The results of the experiments are provided below.
The vertical axis of both  Figure \ref{coord1} and \ref{coord2}  measures function value $f(\rev{y}^k)$ in logarithmic scale, the horizontal axis of Figure \ref{coord1} and \ref{coord2} measures wall-clock time as the algorithm runs. As we see, our algorithm denoted as SAE, performs better when the accuracy is high.
Figure \ref{grad_calc} shows a three-dimensional plot of the function value $f(\rev{y}^k)$ in logarithmic scale vs the number of $\nabla h(\cdot)$ and $\nabla g(\cdot)$ oracle calls \rev{(convergence curve for every method terminates when some common fixed accuracy is reached, so their lengths differ)} and two-dimensional projections of this plot for the $\nabla h(\cdot)$ and $\nabla g(\cdot)$ oracle calls respectively \rev{(curves on 2d projections are clipped at the number of oracle calls that is the least among the presented methods)}.

The experiments reported in this subsection
were performed on a regular laptop with an Intel i5-5250U 1.60 GHz CPU and 8 GB of RAM, under 64-bit macOS operating system, with the Jupyter Notebook browser environment. The version of Python programming language is 3.8.3, the version of the notebook server is 6.1.1, IPython 7.17.0.

Next we also compare theoretical complexity bounds for our approach and other methods used in the experiments.
We assume that in \eqref{eq:softmax_reg} the matrix $A = [A_1,...,A_p]^{\top}$ is such that $\max_{ij} |A_{ij}| = O(1)$, $\max_{j=1,...,n} \|A^{\langle j\rangle}\|_2^2 = O(n)$ and $A$ has $O(ps)$ nonzero elements. Further, we assume that $G^2$ is a positive semidefinite  matrix with $\lambda_{\max}(G^2) = O(n)$ and $\frac{1}{n}\sum_{i=1}^n \sqrt{G^2_{ii}} = O(1)$. We also assume that for the starting point $x_0$ and the closest solution $x^*$ it holds that $\|x_0-x^*\|_2\leq R$.
Under these assumptions, FGM \cite{nesterov2018lectures} requires 
$
O\left(\tfrac{1}{\sqrt{\varepsilon}}\sqrt{\left(\max_{j = 1,..., n} \|A^{\langle j\rangle}\|_2^2 + \lambda_{\max}(G^2)\right)R^2}\right)
$
iterations with the complexity of each iteration being
$O\left(ps + n^2\right)$ arithmetic operations.  ARCD/FCD \cite{nesterov2017efficiency}
requires 
$$
O\left(\frac{n}{\sqrt{\varepsilon}}\sqrt{\left(\max_{ij} |A_{ij}|^2+ \left(\frac{1}{n}\sum_{i=1}^n \sqrt{G^2_{ii}} \right)^2\right)R^2}\right)
$$
iterations with the complexity of each iteration
\footnote{
In this case an efficient way to recalculate the partial derivatives of $h(x))$ is as follows. From the structure of the method we know that $x^{new} = \alpha x^{old} + \beta e_i$, where $e_i$ is $i$-th coordinate vector. Thus, given $\left\langle A_k, x^{old} \right\rangle$ recalculating $\left\langle A_k, x^{new} \right\rangle = \alpha\left\langle A_k, x^{old} \right\rangle + \beta [A_k]_i$ requires only $O(1)$ additional arithmetic operations independently of $n$ and $s$.
}
$O\left(p + n\right)$ arithmetic operations.
For our approach it is sufficient to make 
$
O\left(\tfrac{1}{\sqrt{\varepsilon}}\sqrt{\left(\max_{j = 1,..., n} \|A^{\langle j\rangle}\|_2^2\right)R^2}\right)
$
calls to the full gradient oracle $\rev{\nabla} h$  
with the complexity of each call being
$O(ps)$ arithmetic operations and 
$
\tilde{O}\left(\tfrac{n}{\sqrt{\varepsilon}}\sqrt{ \left(\tfrac{1}{n}\sum_{i=1}^n \sqrt{G^2_{ii}} \right)^2R^2}\right)
$
calls to the coordinate derivative oracle for the quadratic function $g$ with the complexity of each iteration being $O(n)$ arithmetic operations. \rev{The upper bounds for these} results are summarized in the table below. The main observation is that when $n \ll p$ and $s\ll \min\left\lbrace n^2/p,\sqrt{n}\right\rbrace$,  our approach leads to better theoretical complexity. \rev{We notice also that the bounds in the table are upper bounds and may happen to be not tight. This, in particular, means that the gain in practice by our algorithm may be even larger.}
\begin{table}[ht]
\label{table_1}
\vspace{.1 in}
\begin{center}
\begin{tabular}{|c|c|c|ll}
\cline{1-3}
 Algorithm  & Complexity&  Reference&   &  \\ \cline{1-3}
FGM & $ O\left(\sqrt{\frac{nR^2}{\varepsilon}}\left(ps + n^2\right)\right) $ & \cite{nesterov2018lectures} &  &  \\ \cline{1-3}
ARCD/FCD & $ O\left(n\sqrt{\frac{R^2}{\varepsilon}}\left(p + n\right)\right) $ & \cite{nesterov2017efficiency} &  &  \\ \cline{1-3}
\aai{SAE} & $ \tilde{O}\left(\sqrt{\frac{nR^2}{\varepsilon}}\cdot ps\right) + \tilde{O}\left(n\sqrt{\frac{R^2}{\varepsilon}}\cdot n\right) $ & this paper &  &  \\ \cline{1-3}
\end{tabular}
\caption{Comparison of theoretical complexity bounds for \aai{SAE} (Algorithm~\ref{alg:MS}), FGM and ARCD/FCD.}
\end{center}
\end{table}

\section{Proof of the Main Result}\label{App:proof}
In this section we prove Theorem~\ref{main_theo} in the following four steps:
\begin{enumerate}
    \item Estimating the number of iterations of the inner method $\mathcal{M}_{inn}$.
    \item Estimating the number of iterations of the procedure GMCO in Algorithm~\ref{alg:MS}.
    \item Estimating the number of iterations of Algorithm~\ref{alg:MS} on its own and combined with the restarting strategy outlined in Algorithm~\ref{Alg:rest_str_MS}.
    \item Obtaining the final estimates for the number of calls to the oracles $O_f$ and $O_g$ based on the estimates obtained in steps $1-3$.
\end{enumerate}

\textbf{Step 1.}
At each iteration  of the procedure GMCO in Algorithm~\ref{alg:MS} (see steps \ref{step:GMCO_subproblem} and \ref{eq:GMCO_step_1}) we need to solve the following problem (cf. \eqref{eq:innermost_problem}, \eqref{eq:GMCO_step})
\begin{equation}
    \min\limits_{v \in \mathbb{R}^n} \left\{\varphi(v) :=  \left\langle \nabla h(v^{k-1}), v - v^{k-1} \right\rangle 
+ g(v)  + \tfrac{L}{2}\|v-v^0\|_2^2   + \tfrac{L_h}{2} \|v - v^{k-1}\|_2^2\right\}.
    \label{aux_problem}
\end{equation}
By Assumption \ref{ass_inn_method}, applying the method $\mathcal{M}_{inn}$ to~\eqref{aux_problem} with $\alpha = L + L_h$ we obtain that, given some $\tilde{\varepsilon} > 0$, in $N_{\mathcal{M}}(\tilde{\varepsilon}) = O\left(\tfrac{\tau_{g}}{\sqrt{L + L_h}}\ln\tfrac{C \|v^0-v^{\ast}\|_2^2}{\tilde{\varepsilon}} \right)$  calls to the basic oracle $O_g$ we can find $v^{N_{\mathcal{M}}(\tilde{\varepsilon})}$ such that  
\begin{equation}
    \label{varphi_est}
	\E (\varphi(v^{N_{\mathcal{M}}(\tilde{\varepsilon})}) - \varphi(v^{\ast})) \leq \tilde{\varepsilon}.
\end{equation}
Since $\varphi(v^{N_{\mathcal{M}}(\tilde{\varepsilon})}) - \varphi(v^{\ast}) \geq 0$, given an arbitrary $\tilde{\delta}\in(0,1)$ we can apply the Markov inequality and obtain
\begin{equation*}
	\PP \left(\varphi(v^{N_{\mathcal{M}}(\tilde{\delta}\tilde{\varepsilon})}) - \varphi(v^{\ast}) \geq \tilde{\varepsilon} \right) \leq \tfrac{\E (\varphi(v^{N_{\mathcal{M}}(\tilde{\delta} \tilde{\varepsilon})}) - \varphi(v^{\ast}))}{\tilde{\varepsilon}} \leq  \tilde{\delta}.
\end{equation*}
Thus, in 
$N_{\mathcal{M}}(\tilde{\delta}\tilde{\varepsilon}) = O\left(\tfrac{\tau_{g}}{\sqrt{L + L_h}}\ln\tfrac{C \|v^0-v^{\ast}\|_2^2}{\tilde{\delta} \tilde{\varepsilon}} \right)$  calls to the basic oracle $O_g$ we can find $\hat{v}$ such that $\varphi(\hat{v}) - \varphi(v^{\ast}) \leq \tilde{\varepsilon}$ with probability at least $1 - \tilde{\delta}$.
Since $\varphi(\cdot)$ is $L+L_h$ strongly convex, we have
$\tfrac{L+L_h}{2}\|\hat{v} - v^{\ast}\|_2^2 \leq \varphi(\hat{v}) - \varphi(v^{\ast}) \leq \tilde{\varepsilon}$ with probability at least $1 - \tilde{\delta}$.
By the construction of Algorithm \ref{alg:MS} and assumptions of  Theorem \ref{main_theo} we have that $L \leq L_h \leq L_g$, and, hence, $L_g+L+L_h \leq 3 L_g$. Thus, $\varphi(\cdot)$ is $3L_g$-smooth, and since $\nabla \varphi(v^{\ast}) = 0$, we have that  
$\|\nabla \varphi(\hat{v}) \|_2 \leq 3L_g\|\hat{v} - v^{\ast}\|_2$.
Whence, with probability at least $1 - \tilde{\delta}$,
 \begin{eqnarray*}
     \left\langle  \nabla \varphi(\hat{v}), \hat{v} - v^{\ast} \right\rangle  \leq   \|  \nabla \varphi(\hat{v})\|_2 \cdot \| \hat{v} - v^{\ast}  \|_2 \leq 3L_g\|\hat{v} - v^{\ast}\|_2^2 \leq 3L_g \tfrac{2\tilde{\varepsilon}}{L+L_h}.
 \end{eqnarray*}
 This leads to the following lemma that summarizes Step 1.
 \begin{lemma}
 \label{inn_est}
 Applying $\mathcal{M}_{inn}$ to~\eqref{aux_problem} we have that  in \\
 $N_{\mathcal{M}}(\tilde{\delta}\tilde{\varepsilon}) = O\left(\tfrac{\tau_{g}}{\sqrt{L + L_h}}\ln\tfrac{2C L_g \|v^0-v^{\ast}\|_2^2}{\tilde{\delta}\tilde{\varepsilon}(L+L_h)} \right)$ calls to the basic oracle $O_g$ we can find $\hat{v}$ s.t. with probability at least $1 - \tilde{\delta}$
 \begin{align}
  \label{inex_sol}
    \left\langle  \nabla \varphi(\hat{v}), \hat{v} - v^{\ast} \right\rangle \leq \tilde{\varepsilon}, \qquad \|\hat{v} - v^{\ast}\|^2_2 \leq \tfrac{\tilde{\varepsilon}}{L_g}.
    \end{align}
 \end{lemma}

\textbf{Step 2.}
Let $\zeta^{\ast}$ be the exact minimizer of the function  $F_{L, \zeta^0}(\zeta) = f(\zeta) + \tfrac{L}{2} \|\zeta - \zeta^0\|^2_2$ that is minimized by the procedure GMCO in Algorithm~\ref{alg:MS} (see steps \ref{step:applyGMCO} and \ref{step:GMCO_input}). 
Let also $\hat{\zeta}^k$ be an inexact minimizer of this function that satisfies the 
Monteiro-Svaiter (MS) condition
 $ \|\nabla F_{L, \zeta^0 } (\hat{\zeta}^k)\|_2 \leq \tfrac{L}{2} \|\hat{\zeta}^k - \zeta^0\|_2 $ (cf. step \ref{step:checkMS} of GMCO).
By the optimality condition, $  \nabla F_{L, \zeta^0}  (\zeta^{\ast})  = 0$, and, thus,  
\begin{equation}
\label{eq:1}
    \| \nabla F_{L, \zeta^0}(\hat{\zeta}^k)\|_2 \leq  (L+L_f)\|\hat{\zeta}^k - \zeta^{\ast}\|_2,
\end{equation}
since  $F_{L, \zeta^0}(\cdot)$ is $(L+L_f)$-smooth.
By the triangle inequality we get 
\begin{equation}
  \label{eq:2}
    \frac{L}{2}\|\zeta^0 - \zeta^{\ast}\|_2 - \frac{L}{2}\|\hat{\zeta}^k -  \zeta^{\ast}\|_2 \leq \frac{L}{2}\| \hat{\zeta}^k - \zeta^0 \|_2.  
\end{equation}
Since the r.h.s. of \eqref{eq:2} coincides with the r.h.s. of the MS condition and l.h.s. of \eqref{eq:1}
coincides with the l.h.s. of the MS condition, we see that if the inequality 
\begin{equation}
\label{inner_est}
      \|\hat{\zeta}^k - \zeta^{\ast}\|_2 \leq \tfrac{L}{3 L + 2L_{f}}  \|\zeta^0 -\zeta^{\ast}\|_2
\end{equation}
holds, the MS condition holds too. Thus, our next goal is to ensure condition \eqref{inner_est} by the procedure GMCO in Algorithm \ref{alg:MS}. This in turn will imply that MS condition in step \ref{step:checkMS} of the procedure GMCO holds.

To that end, we assume that at each iteration of the procedure GMCO we solve the auxiliary problem of minimizing the function~\eqref{eq:GMCO_step} (see also \eqref{aux_problem}) in the sense of the left inequality in \eqref{inex_sol}. 
\aai{Under this assumption, the following convergence rate theorem  for the procedure GMCO is a direct extension of Theorem  3.1 in \cite{stonyakin2020inexact}. The full proof of this theorem is provided in Section \ref{GMCO}.}

\begin{theorem}\label{Th:MS}
Let $\frac{\mu + L}{2L_h} \leq 1$. Then, under the above assumption on inexact minimization of \eqref{eq:GMCO_step}, after $K$ iterations of the procedure GMCO in Algorithm \ref{alg:MS} we have
\begin{eqnarray}
 F_{L, \zeta^0}(\zeta^K) - F_{L, \zeta^0}(\zeta^{\ast}) &\leq& \exp\left(-\tfrac{K(\mu + L)}{4L_h}\right) (F_{L, \zeta^0}(\zeta^0) - F_{L, \zeta^0}(\zeta^{\ast})) + \tfrac{4L_h}{\mu + L}\tilde{\varepsilon}, \notag \\
 \tfrac{1}{2}\|\zeta^{\ast} - \zeta^K\|^2_2 &\leq& \tfrac{L_h}{2(\mu + L)} \|\zeta^{\ast} -  \zeta^0\|^2_2 + \tfrac{4L_h}{(\mu + L)^2}\tilde{\varepsilon}.
\label{eq_zeta_est}
\end{eqnarray}
\end{theorem}


Since  $F_{L, \zeta^0}(\cdot)$ is $L$-strongly convex, by Theorem~\ref{Th:MS}, we have that 
\begin{equation*}
    \tfrac{L}{2}\| \zeta^K- \zeta^{\ast}\|^2_2 \leq F_{L, \zeta^0}(\zeta^K) - F_{L, \zeta^0}(\zeta^{\ast} )   \leq \tfrac{L_{h}\|\zeta^0 - \zeta^{\ast}\|_2^2}{2}\exp\left( \tfrac{-K L } {4L_h}\right) + \tfrac{4L_h}{L}\tilde{\varepsilon}.
\end{equation*}
Thus, for the condition~\eqref{inner_est} to be satisfied by the point $ \zeta^K$ produced by the procedure GMCO, it is sufficient that
\begin{equation*}
     \tfrac{L_{h}\|\zeta^0 - \zeta^{\ast}\|_2^2}{2}\exp\left( \tfrac{-K L } {L_h}\right) + \tfrac{4L_h}{L}\tilde{\varepsilon} \leq \tfrac{L^3}{2(3 L + 2L_{f})^2} \|\zeta^0 - \zeta^{\ast}\|^2_2.
\end{equation*}
Equating each term of the l.h.s. to half of the r.h.s., we obtain that a sufficient number of iterations of the the procedure GMCO is

\noindent$
    K_{\text{GMCO}} := \tfrac{L_h}{L} \ln\left(\tfrac{2(3 L + 2L_{f})^2 L_h}{L^3} \right) 
$
and a sufficient accuracy $\tilde{\varepsilon}$ for the solution of the auxiliary problem of minimizing the function ~\eqref{eq:GMCO_step} (see also \eqref{aux_problem}) in the sense of left inequality in \eqref{inex_sol} is 
\begin{equation}
    \label{eq:eps_M_inn}
    \tilde{\varepsilon} = \varepsilon_{\mathcal{M}} := \tfrac{L^4}{16L_h(3 L + 2L_{f})^2} \|\zeta^0 - \zeta^{\ast}\|^2_2.
\end{equation}
Further, let $\delta_{\text{SAE}} \in (0,1)$ be some confidence probability level and assume that in each iteration of the procedure GMCO an $\varepsilon_{\mathcal{M}}$-solution to the auxiliary problem of minimizing function \eqref{eq:GMCO_step} is obtained with probability at least   $1 - \delta_{\text{SAE}} / K_{\text{GMCO}}$. Then, using the union bound, Step 2 can be summarized by the following result.
\begin{lemma}
\label{Lm:GMCO_compl}
   In $K_{\text{GMCO}}$ iterations of the the procedure GMCO in Algorithm~\ref{alg:MS}, assuming that in each iteration an $\varepsilon_{\mathcal{M}}$-solution to the auxiliary problem is obtained with probability at least   $1 - \delta_{\text{SAE}} / K_{\text{GMCO}}$,    we find a point $\hat{\zeta}$ such that
   $    \|\nabla F_{L, \zeta^0 } (\hat{\zeta})\|_2 \leq \tfrac{L}{2} \|\hat{\zeta} - \zeta^0\|_2$ with probability at least $1 - \delta_{\text{SAE}}$.
\end{lemma}

\textbf{Step 3.}
To estimate the number of iterations of SAE (Algorithm~\ref{alg:MS}) note that in~\eqref{aux_pr} we apply the procedure GMCO and, according to the stopping criterion of  $\text{GMCO}(x^{k+1},F_{L, x^{k+1}})$ in step \ref{step:checkMS} obtain $y^{k+1}$ such that 
\begin{equation}
\label{MS_cond}
  \|\nabla F_{L, x^{k+1}} (y^{k+1})\|_2 \leq \tfrac{L}{2} \|y^{k+1} - x^{k+1}\|_2.   
\end{equation}
Applying  Theorem 3.6 from~\cite{monteiro2013accelerated} to  SAE, we obtain that, for all $N \geq 0$,
\begin{equation}
\label{eq12}
f\left( {y^N} \right)-f\left( {x^{\ast} } \right)\le \tfrac{\|x^0 - x^{\ast}\|_2^2}{2A_N }\leq \tfrac{R^2}{2A_N },
\quad
\left\| {z^N-x^{\ast} } \right\|_2 \le R,
\end{equation}
where $R \geq \left\| {y^0-x^{\ast} } \right\|_2$. Moreover, from Lemma 3.7 a) of ~\cite{monteiro2013accelerated}, for all $N \geq 0$, $A_N \ge \tfrac{N^2}{4L}$.
Substituting the last inequality into the estimate \eqref{eq12}, we obtain that after $N$ iterations of SAE the following inequality holds
 \begin{equation}
 \label{MS_func_est}
     f(y^N) - f(x_\ast) \leq \tfrac{2L\|x^0 - x^{\ast}\|_2^2}{N^2}.
 \end{equation}
 Thus,  if $\mu = 0$, then the total number of SAE iterations  to achieve accuracy $\varepsilon$ is $T^{\text{c}}_{\text{SAE}}(\varepsilon) := \sqrt{\tfrac{LR^2}{\varepsilon}}$.
 
 If $\mu > 0$ we need to apply the restarting strategy described in Algorithm \ref{Alg:rest_str_MS}. 
By the strong convexity of $f(\cdot)$ and \eqref{MS_func_est}, we have
 \begin{eqnarray*}
    \tfrac{\mu}{2}\|y^N - x^{\ast}\|_2 \leq  f(y^N) - f(x^{\ast}) \leq \tfrac{2L}{N^2}\|x^0 - x^{\ast}\|_2^2.
 \end{eqnarray*}
 Thus, after $N = N_0 = \sqrt{\tfrac{8L}{\mu}}$ iterations, we obtain $\|y^N - x^{\ast}\|_2 \leq \tfrac{1}{2}\|x^0 - x^{\ast}\|_2^2$ and restart Algorithm \ref{alg:MS}. Hence, after $t$ restarts we have
  \begin{eqnarray*}
     f(\eta_{t+1}) - f(x^{\ast}) = f(y^{N_0}) - f(x^{\ast}) \leq \tfrac{4 L}{2^t N_{0}^2}\|x^0 - x^{\ast}\|_2^2 = \tfrac{\mu}{2^{t+1}}\|x^0 - x^{\ast}\|_2^2.
 \end{eqnarray*}
By choosing $t=T \geq \ln\left(\tfrac{\mu\|x_0 - x^{\ast}\|_2^2 }{\varepsilon} \right)$, we see that after $T$
iterations of  Algorithm~\ref{Alg:rest_str_MS} we obtain $\eta_T$ such that $f(\eta_T) - f(\eta^{\ast}) \leq \varepsilon$.
Since each iteration of Algorithm~\ref{Alg:rest_str_MS} uses $N_0=\sqrt{\tfrac{8L}{\mu}}$ iterations of SAE, we obtain that if $\mu >0$, the total number of $\text{SAE}$ iterations to find an $\varepsilon$-solution to problem \eqref{prob_st} is $T^{\text{sc}}_{\text{SAE}}(\varepsilon):= O\left(\sqrt{\tfrac{L}{\mu}} \ln\left(\tfrac{\mu\|\eta_0 - \eta^{\ast}\|_2^2}{\varepsilon} \right)\right)$. 

Assume now that at each iteration of Algorithm~\ref{alg:MS} we find $y^{k+1}$ satisfying~\eqref{MS_cond} with probability at least $1- \delta / T_{\text{SAE}}(\varepsilon)$, where $\delta \in (0,1)$ and
$$
T_{\text{SAE}}(\varepsilon) = 
  \begin{cases} 
   T^{\text{sc}}_{\text{SAE}}(\varepsilon) & \text{if } \mu > 0, \\
   T^{\text{c}}_{\text{SAE}}(\varepsilon)       & \text{if } \mu = 0.
  \end{cases}
 $$ 
  
Using the union bound over all $\text{SAE}$ iterations we summarize this step in the following lemma.
\begin{lemma}
\label{Lm:SAE_compl}
    If in each iteration of Algorithm~\ref{alg:MS} we find $y^{k+1}$ satisfying~\eqref{MS_cond} with probability at least $1-\delta_{\text{SAE}}$ with $\delta_{\text{SAE}} = \delta / T_{\text{SAE}}(\varepsilon)$, then
    
    a) after $T^{\text{sc}}_{\text{SAE}}(\varepsilon)$ iterations of Algorithm~\ref{alg:MS} for the case $\mu > 0$ 
    
    b) after $T^{\text{c}}_{\text{SAE}}(\varepsilon)$ iterations of  Algorithm~\ref{alg:MS} for the case $\mu = 0$ 
    
    we find a point $\hat{\eta}$ such that $f(\hat{\eta}) - f(\eta^{\ast}) \leq \varepsilon$ with probability as least $1 - \delta$.
\end{lemma}

\textbf{Step 4.} 
At each iteration of  the procedure GMCO in Algorithm \ref{alg:MS} we use the method $\mathcal{M}_{inn}$ with starting point $\zeta^k$ to compute the point $\zeta^{k+1}$. So, for the $k$-th iteration of  GMCO, we apply Lemma \ref{inn_est} with $v^0 \equiv \zeta^k$ and $\hat{v} \equiv \zeta^{k+1}$.
Then, using the triangle inequality and \eqref{eq_zeta_est}, we get
\begin{eqnarray*}
    \|\zeta^k - \zeta^{k+1}\|_2 &\leq& \|\zeta^k - \zeta^{\ast}\|_2 + \|\zeta^{k+1} - \zeta^{\ast}\|_2 \\ &\leq& 2\|\zeta^k - \zeta^{\ast}\|_2  \leq 2 \sqrt{\tfrac{L_h}{L}} \|\zeta^0 - \zeta^{\ast}\|_2 + 4 \sqrt{\tfrac{2L_h}{L^2}\tilde{\varepsilon}}.
\end{eqnarray*}
From this, by using the right inequality in~\eqref{inex_sol}, we have 
\begin{equation}
      \| v^0 - v^{\ast}\|_2  \leq  \|v^0 - \hat{v}\|_2 + \sqrt{\frac{\tilde{\varepsilon}}{L_g}  }   \leq  \sqrt{\tfrac{4L_h}{L}} \|\zeta^0 - \zeta^{\ast}\|_2 + \sqrt{\tilde{\varepsilon}} \left( \sqrt{\tfrac{32L_h}{L^2}} + \sqrt{\tfrac{1}{L_g}  }  \right).
     \label{norm_est}
\end{equation}
Choosing $\tilde{\varepsilon} = \varepsilon_{\mathcal{M}}$ defined in \eqref{eq:eps_M_inn} and using Lemma~\ref{inn_est}, we obtain that we need   $N_{\mathcal{M}}(\varepsilon_{\mathcal{M}}) = O\left(\tfrac{\tau_{g}}{\sqrt{L + L_h}}\ln\tfrac{2C L_g \|v^0 - v^{\ast}\|_2^2 }{\tilde{\delta}\varepsilon_{\mathcal{M}} (L+L_h)} \right)$ calls to the oracle $O_g$ used by $\mathcal{M}_{inn}$ to find $\hat{v}$ s.t. with probability at least $1 - \tilde{\delta}$  conditions \eqref{inex_sol} hold.
Combining the definition \eqref{eq:eps_M_inn} of $\varepsilon_{\mathcal{M}}$ with \eqref{norm_est} and denoting $C_1 = \tfrac{2C(32L_h^2 L_g (3L+2L_f)^2 + (8L_h L_g + L^2)L^3}{L^4 (L+L_h)}$, we obtain the following upper bound  $\tilde{N}_{\mathcal{M}} = O\left(\tfrac{\tau_{g}}{\sqrt{L + L_h}} \ln\tfrac{C_1  }{\tilde{\delta}}  \right) \geq N_{\mathcal{M}}(\varepsilon_{\mathcal{M}})$. 
   
Based on Lemmas \ref{Lm:GMCO_compl} and \ref{Lm:SAE_compl}, we  choose $\tilde{\delta} \approx  \delta_{\text{SAE}}/ K_{\text{GMCO}} = \delta / (K_{\text{GMCO}} \cdot T_{\text{SAE}})$ and obtain that, in the innermost cycle, to find an approximate minimizer of \eqref{eq:GMCO_step} (cf. \eqref{aux_problem}) $\hat{v}$ s.t. with probability at least $1 - \tilde{\delta}$ inequalities \eqref{inex_sol} hold, we need $N^{\text{sc}}_{\mathcal{M}} = O\left(\tfrac{\tau_{g}}{\sqrt{L + L_h}} \ln\tfrac{C_1 L_h }{\delta \sqrt{\mu L}}  \right)$  calls to the oracle $O_g$ used by $\mathcal{M}_{inn}$ for the strongly convex case ($\mu >0$) and $N^{\text{c}}_{\mathcal{M}} = O\left(\tfrac{\tau_{g}}{\sqrt{L + L_h}} \ln\tfrac{C_1 L_h R}{\delta \sqrt{\varepsilon L}}  \right)$ calls to the oracle $O_g$ used by $\mathcal{M}_{inn}$ for the convex case ($\mu=0$), where $R \geq \|x^0 - x^{\ast}\|_2$.
   
To obtain the estimates for the total number of oracle calls for  $h(\cdot)$ and $g(\cdot)$, we make the following observations. 
We need to compute the gradient of  $h(\cdot)$ at each step of the procedure GMCO in Algorithm~\ref{alg:MS}, that we run in total $T_{\text{SAE}}(\varepsilon)$ times. Moreover, at each iteration of  Algorithm~\ref{alg:MS} in step 5 we compute the gradient of $f(\cdot)$ that contains the gradient of $h(\cdot)$.
Further,  we need to call the basic oracle for $g(\cdot)$ at each step of the inner algorithm $\mathcal{M}_{inn}$, that we run at each iteration of the procedure GMCO, and at each iteration of  Algorithm~\ref{alg:MS} in step 5 we need to calculate the full gradient of $g(\cdot)$ that according to Assumption \ref{Asmpt:g} costs $\kappa_{g}$ calls to the basic oracle $O_g$.   
   
   Using the union bound over all launches of $\mathcal{M}_{inn}$, we obtain that our Algorithms \ref{alg:MS} and \ref{Alg:rest_str_MS} find such $\hat{x}$ that $f(\hat{x}) - f(x^{\ast}) \leq \varepsilon$ with probability at least $1 - \delta$, and to do this it is sufficient to make $O(T_{\text{SAE}} \cdot (1+K_{\text{GMCO}}))$ oracle calls for $h(\cdot)$ and $O\left(T_{\text{SAE}} \cdot \left(\kappa_g + K_{\text{GMCO}} \cdot  N_{\mathcal{M}}\right)\right)$ oracle calls for $g(\cdot)$.
   
   So, we need 
   $
   O\left(\sqrt{\tfrac{L}{\mu }} \ln\left(\tfrac{  \mu R^2}{\varepsilon} \right) \cdot \left(1 + \tfrac{L_h}{L} \right) \right)
 $
oracle calls for $h(\cdot)$ and \begin{eqnarray*}
   O \left(\sqrt{\tfrac{L}{\mu}}  \ln\left(\tfrac{  \mu R^2}{\varepsilon}\right) \cdot \left(\kappa_g + \tfrac{L_h}{L}\cdot \left(\tfrac{\tau_{g}}{\sqrt{L + L_h}}\ln\tfrac{1}{\delta } \right) \right) \right) 
   \end{eqnarray*}
   oracle calls for $g(\cdot)$,  if $\mu > 0$,
   
   and 
$
   O\left(\sqrt{\tfrac{LR^2}{\varepsilon}} \cdot \left(1 + \tfrac{L_h}{L} \right) \right)$
oracle calls for $h(\cdot)$ and 
\begin{eqnarray*}
   O \left(\sqrt{\tfrac{LR^2}{\varepsilon}} \cdot \left(\kappa_g + \tfrac{L_h}{L}\cdot \left(\tfrac{\tau_{g}}{\sqrt{L + L_h}}\ln\tfrac{R}{\delta \sqrt{\varepsilon}} \right) \right) \right) 
   \end{eqnarray*}
   oracle calls for $g(\cdot)$, if $\mu = 0$.
This finishes the proof of Theorem \ref{main_theo}.

\section{\rev{Proof of Theorem \ref{Th:MS}}}\label{GMCO}
\rev{
We start by a definition of inexact solution to the auxiliary problem in Gradient Method for Composite Optimization.
\begin{definition}
\label{solNemirovskiy}
For a convex optimization problem 
$\min_{x \in Q} \Psi(x)$,
we denote by $\text{Arg}\min_{x \in Q}^{\widetilde{\delta}}\Psi(x)$~a set of such $\widetilde{x}$ that
\begin{gather}\label{eqv_inex_sol}
\exists h \in \partial\Psi(\widetilde{x}):\forall x \in Q \,\, \to\, \left\langle h, x - \widetilde{x} \right\rangle \geq -\widetilde{\delta}.
\end{gather}
We denote by $\argmin_{x \in Q}^{\widetilde{\delta}}\Psi(x)$  some element of $\text{Arg}\min_{x \in Q}^{\widetilde{\delta}}\Psi(x)$.
\end{definition}
\begin{algorithm}[ht]
	\caption{Gradient Method for Composite Optimization $\text{GMCO}(x_0, F(\cdot))$}
	\begin{algorithmic}[1]
		\STATE {\bf Parameters:} starting point $x_0 \in \mathbb{R}^n$, objective function $F(x) = f(x) + p(x)$, constant $L$ (function $f$ with $L$ Lipschitz gradient w.r.t. the $||\cdot||_2$), error $\widetilde{\delta}$.
		\FOR{$k = 0, \ldots,N-1$}
			\STATE Set
			$
			\phi_{k+1}(x) := \left\langle \nabla f(x_k), x - x_k \right\rangle + p(x) + \tfrac{L}{2} \|x - x_k\|_2^2,
			$
			\STATE Compute 
			\begin{equation}
			\label{alg:gradient_for_composite_base:auxilary}
			   x_{k+1}:={\argmin_{x \in Q}}^{\widetilde{\delta}}(\phi_{k+1}(x))
			\end{equation}
		\ENDFOR
		\STATE {\bf Output:} $x_N$
	\end{algorithmic}
	\label{alg:gradient_for_composite_base}
\end{algorithm}
\begin{lemma}
    \label{lemma:str_conv}
	Let $\psi(x)$ be a convex function and
	\begin{gather*}
	y = {\argmin_{x \in Q}}^{\widetilde{\delta}} \left\lbrace \psi(x) + \tfrac{\beta}{2}||z - x||^2_2\right\rbrace,
	\end{gather*}
	where $\beta \geq 0$.
Then
	\begin{align*}
\psi(x) + \tfrac{\beta}{2}||z - x||^2_2
\geq \psi(y) + \tfrac{\beta}{2}||z - y||^2_2 + \tfrac{\beta}{2}||x - y||^2_2 - \widetilde{\delta} ,\,\,\, \forall x \in Q.
	\end{align*}
\end{lemma}
\begin{proof}
    By Definition \ref{solNemirovskiy}:
    \begin{gather*}
		\exists g \in \partial\psi(y), \,\,\, \left\langle g + \tfrac{\beta}{2} \nabla_y ||y - x||^2_2, x - y \right\rangle = \left\langle g + \beta(y - z), x - y \right\rangle \geq -\widetilde{\delta} ,\,\,\, \forall x \in Q.
	\end{gather*}
    From $\beta$--strong convexity of $\psi(x) + \tfrac{\beta}{2}||z - x||^2_2$ we have
    \begin{gather*}
		\psi(x) + \tfrac{\beta}{2}||z - x||^2_2 \geq \psi(y) + \tfrac{\beta}{2}||z - y||^2_2 + \left\langle g + \tfrac{\beta}{2} \nabla_y ||y - x||^2_2, x - y\right\rangle + \tfrac{\beta}{2}||x - y||^2_2
	\end{gather*}
	The last two inequalities complete the proof.
\end{proof}
The next theorem proves convergence rate of Algorithm \ref{alg:gradient_for_composite_base} for optimization problem
\begin{equation}
    \min_{x \in \mathbb{R}^n} F(x):= f(x) + p(x),
\end{equation}
where function $f$ is convex function with $L$-Lipschitz gradient w.r.t. the $||\cdot||_2$ norm, function $p$ is convex function and function $F$ is $\mu$--strongly convex.
\begin{theorem}\label{Th:str_conv_adap}
Let us assume that $\frac{\mu}{2L} \leq 1$. After $N$ iterations of Algorithm \ref{alg:gradient_for_composite_base} we have
\begin{align*}
F(x_{N}) - F(x^{\ast}) \leq \exp\left(-\frac{N\mu}{4L}\right) (F(x_{0}) - F(x^{\ast})) + \frac{4L}{\mu}\widetilde{\delta},
\end{align*}
\begin{align*}
    \frac{1}{2}||x^{\ast} - x_{N}||^2_2 \leq \frac{L}{2\mu} ||x^{\ast} - x_{0}||^2_2 + \frac{4L}{\mu^2}\widetilde{\delta}.
\end{align*}
\end{theorem}
{\it Proof}
Since gradient of function $F$ is $L$-Lipschitz w.r.t. the $||\cdot||_2$ norm, we have
\begin{align*}
F(x_{N}) \leq  f(x_{N-1}) + \left\langle \nabla f(x_{N-1}), x_{N} - x_{N-1} \right\rangle + p(x_N) + \frac{L}{2}||x_{N-1} - x_{N}||^2_2.
\end{align*}
From Lemma \ref{lemma:str_conv} and auxiliary problem \eqref{alg:gradient_for_composite_base:auxilary} we get
\begin{align*}
F(x_{N})
\leq f(x_{N-1}) +\left\langle \nabla f(x_{N-1}), x - x_{N-1} \right\rangle + p(x) + \frac{L}{2}||x - x_{N-1}||^2_2 + \widetilde{\delta}.
\end{align*}
In view of convexity of function $f$, we obtain
\begin{align}
\label{analysis_algorithm_grad:prove_1}
F(x_{N}) \leq F(x) + \frac{L}{2}||x - x_{N-1}||^2_2 + \widetilde{\delta}.
\end{align}
We rewrite the last inequality for $x = \alpha x^{\ast} + (1 - \alpha) x_{N-1}$ ($\alpha \in [0, 1]$) as
\begin{align*}
F(x_{N}) \leq F( \alpha x^{\ast} + (1 - \alpha) x_{N-1}) + \frac{L\alpha^2}{2}||x^{\ast} - x_{N-1}||^2_2 + \widetilde{\delta}.
\end{align*}
In view of convexity of function $f$, we have
\begin{align*}
F(x_{N}) \leq F(x_{N-1}) -\alpha (F(x_{N-1}) - F(x^{\ast})) + \frac{L\alpha^2}{2}||x^{\ast} - x_{N-1}||^2_2 + \widetilde{\delta}.
\end{align*}
From $\mu$--strong convexity of function $F$ we have $F(x_{N-1}) \geq F(x^{\ast}) + \frac{\mu}{2}||x^{\ast} - x_{N-1}||^2_2$, this yields inequality:
\begin{align*}
F(x_{N}) \leq F(x_{N-1}) -\alpha \left(1 - \alpha\frac{L}{\mu}\right) (F(x_{N-1}) - F(x^{\ast})) + \widetilde{\delta}.
\end{align*}
The minimum of the right part of the last inequality is achieved with $\alpha = \min(1, \frac{\mu}{2L})$. Due to $\frac{\mu}{2L} \leq 1$ with $\alpha = \frac{\mu}{2L}$ we have
\begin{align*}
F(x_{N}) - F(x^{\ast}) \leq \left(1 - \frac{\mu}{4L}\right) (F(x_{N-1}) - F(x^{\ast})) + \widetilde{\delta}.
\end{align*}
and
\begin{align*}
F(x_{N}) - F(x^{\ast}) &\leq \left(1 - \frac{\mu}{4L}\right)^N (F(x_{0}) - F(x^{\ast})) + \frac{4L}{\mu}\widetilde{\delta}\\
&\leq \exp\left(-\frac{N\mu}{4L}\right) (F(x_{0}) - F(x^{\ast})) + \frac{4L}{\mu}\widetilde{\delta}.
\end{align*}
From $\mu$--strong convexity of function $F$ and the fact that gradient of function $F$ is $L$ Lipschitz we obtain
\begin{align*}
    \frac{1}{2}||x^{\ast} - x_{N}||^2_2 &\leq \frac{L}{2\mu} \exp\left(-\frac{N\mu}{4L}\right) ||x^{\ast} - x_{0}||^2_2 + \frac{4L}{\mu^2}\widetilde{\delta}\\
    &\leq \frac{L}{2\mu} ||x^{\ast} - x_{0}||^2_2 + \frac{4L}{\mu^2}\widetilde{\delta}.
\end{align*}
\qed
}
\section{Conclusions}\label{S:concl}
In this paper we consider minimization problems with objective that is expressed as a sum of two terms. We provide a general algorithmic framework for such problems and obtain complexity bounds corresponding to accelerated methods. One of the key features of our framework is separation of the complexities for each component in the sum. We consider several particular cases where our framework can be applied: full-gradient methods, random coordinate descent and variance reduced methods. We illustrate our general framework by two particular experiments in application to Kernel SVM and log-density estimation. We further explain the regime in which our algorithm has better complexity bounds than the existing algorithms for these two applications. In the future it would be interesting to extend our framework to other randomized methods such as random gradient-free methods. 

Since this paper first appeared as a preprint \cite{ivanova2020oracle}, its main results were extended and applied in several research directions. 
One important application is decentralized distributed optimization (see \cite{dvinskikh2021decentralized,rogozin2021accelerated} and references therein), where the function $g$ is the penalty for the consensus constraints. In this setup separating oracle complexities play an important role since it allows to save the number of communication rounds.  
Another line of applications is related to composite saddle-point problems \cite{gasnikov2021accelerated,vladislav2021accelerated} and smooth variational inequalities \cite{rogozin2021decentralized,lan2021mirror}. 
Finally, it seems that such techniques can have a valuable impact not only in application to smooth problems, but also for higher-order smooth, non-smooth and stochastic problems with different orders of the available oracle \cite{kamzolov2020optimal,beznosikov2020derivative,gladin2021solving,stepanov2021one}. So, we see this paper as a paper providing certain perspective for obtaining different classes of optimization algorithms in general.

As a particular potential application we would like to mention robustifying Artificial Intelligence models by studying adversarial attacks. One of the approaches is based on zeroth-order algorithms that are used in the research of robustness of deep neural networks (DNN) in order to construct adversarial examples that are misclassfied by a DNN \cite{chen2017ZOO,tu2019autoZOOM}.
In particular, the authors of \cite{chen2017ZOO} consider the problem of the form 
\begin{equation}
\label{eq:attack}
    \min_{x \in \R^n} h(x) + g(x),
\end{equation}
where the variable $x$ corresponds to the constructed adversarial attack, function $h$ is a smooth function with gradient oracle which is responsible for reduction of the dimension in the attack $x$, function $g$ is defined using the output of a DNN and is available through zeroth-order oracle. 
Such problem structure fits quite well the problem template \eqref{prob_st} and motivates to use some zeroth-order methods \cite{bogolubsky2016learning,nesterov2017random,dvurechensky2017randomized,dvurechensky2018accelerated,vorontsova2019accelerated,dvurechensky2021accelerated,shibaev2021zeroth-order} as $\mathcal{M}_{inn}$. Then, oracle complexity separation will allow to save a large number of evaluations of $\nabla h(x)$ since, typically, theoretical bounds on the number of iterations of existing zeroth-order methods applied to \eqref{eq:attack} are proportional to $n$.

\begin{acknowledgements}
This work was supported by a grant for research centers in the field of artificial intelligence, provided by the Analytical Center for the Government of the Russian Federation in accordance with the subsidy agreement (agreement identifier 000000D730321P5Q0002 ) and the agreement with the Ivannikov Institute for System Programming of the Russian Academy of Sciences dated November 2, 2021 No. 70-2021-00142.
\end{acknowledgements}

%
%


\bibliographystyle{spmpsci}
\bibliography{biblio}

%
%

\end{document}